%% file: low-rank_convexification.tex
\documentclass[11pt,reqno]{amsart}

\input{isomath}

\usepackage{algorithm,algorithmicx}
\usepackage{algpseudocode}
\usepackage{bbm}
\usepackage[mathlines,pagewise,left]{lineno}
\usepackage{amssymb}
\usepackage{color,xcolor,tikz}
\usetikzlibrary{arrows,shapes,chains}  
\usepackage{amsmath}
\usepackage{bcol}
\usepackage{multirow, array}
\usepackage{caption}
\usepackage{xspace}
\usepackage[numbers]{natbib}
\usepackage{comment}
\usepackage[toc,page]{appendix}
\usepackage{graphicx} 
\usepackage{booktabs} 
\usepackage{pdflscape}
\usepackage[export]{adjustbox}
\usepackage[font=scriptsize]{subfig}
\usepackage{mathtools}
\usepackage{hyperref}
\usepackage{amsthm}
\theoremstyle{definition}
\newtheorem{assumption}{Assumption}
\theoremstyle{remark}
\input{notations}

\newcommand{\set}[1]{\ensuremath{\mathcal{#1}}}
\renewcommand{\conv}{\ensuremath{\text{cl conv}}}

\newcommand{\basic}{\ensuremath{\mathsf{Basic}}}
\newcommand{\rankOne}{\ensuremath{\mathsf{RankOne}}}
\newcommand{\rankTwo}{\ensuremath{\mathsf{RankTwo}}}

\newcolumntype{\resetRow}{>{\global\let\currentrowstyle\relax}}
\newcolumntype{^}{>{\currentrowstyle}}

\def\SingleSpacedXI{\linespread{1.1}}
\SingleSpacedXI

\usepackage[colorinlistoftodos,prependcaption,textsize=small]{todonotes}
\title[Low-rank convexification]{Compact extended formulations for low-rank functions with indicator variables}
\author{Shaoning Han and Andr\'es G\'omez}
\thanks{ \noindent \hskip -5mm
	 S. Han, A. G\'{o}mez: Daniel J. Epstein Department of Industrial and Systems Engineering, Viterbi School of Engineering, University of Southern California, CA 90089. \texttt{shaoning@usc.edu}, \texttt{gomezand@usc.edu}. 
}
\begin{document}
\maketitle
\begin{abstract}
	%\vskip 3mm
	\noindent We study the mixed-integer epigraph of a special class of convex functions with non-convex indicator constraints, which are often used to impose logical constraints on the support of the solutions. The class of functions we consider are defined as compositions of low-dimensional nonlinear functions with affine functions Extended formulations describing the convex hull of such sets can easily be constructed via disjunctive programming, although a direct application of this method often yields prohibitively large formulations, whose size is exponential in the number of variables. In this paper, we propose a new disjunctive representation of the sets under study, which leads to compact formulations with size exponential in the dimension of the nonlinear function, but polynomial in the number of variables. Moreover, we show how to project out the additional variables for the case of dimension one, recovering or generalizing known results for the convex hulls of such sets (in the original space of variables). Our computational results indicate that the proposed approach can significantly improve the performance of solvers in structured problems.
	
	\noindent
	\textbf{Keywords}. Mixed-integer nonlinear optimization, convexification, disjunctive programming, indicator variables. \\
\end{abstract}

\section{Introduction}\label{sec:intro}
In this paper, we consider a general mixed-integer convex optimization problem with indicator variables
\begin{equation}\label{eq:minlp}
	\min_{x,z}\;\left\{ F(x):(x,z)\in \set F\subseteq \R^n\times\{0,1\}^n, x_i(1-z_i)=0\;\forall i\in[n] \right\},
\end{equation} where $F:\R^n\to\R$ is a convex function, $\set F$ is the feasible region, and $[n]\defeq\{1,2,\dots,n\}$.  Each binary variable $z_i$ indicates whether a given continuous variable $x_i$ is zero or not. In other words, $z_i=0\implies x_i=0$, and $z_i=1$ allows $x_i$ to take any value. Optimization problem \eqref{eq:minlp} arises in a variety of settings, including best subset selection problems in statistics \cite{bertsimas2016best}, and portfolio optimization problems in finance \cite{bienstock1996computational}.
In practice, the objective function often has the form $F(x)=\sum_{i\in[m]}f_i(x)$, where each $f_i(x)$ is a composition of a relatively simple convex function and a low-rank linear map. {In the context of machine learning and statistics, $m$ represents the number of samples, and in turn, each index $i$ corresponds to one observed data point.  Besides linear regression for which each $f_i$ is a rank-one quadratic function, another prominent class of statistical models is logistic regression --widely used for classification problems--, where each $f_i$ is a composition of a logistic function and a linear function, that is, $f_i(x)=\ln\left(e^{a_i^\top x}+1\right)$ for a suitable vector $a_i\in \R^n$.  Moreover, in some scenarios \cite{rudin2018optimized}, decision variables are additionally required to be nonnegative for better interpretability and possess certain combinatorial structures such as sparsity. Such optimization problems can be readily modeled as \eqref{eq:minlp} \cite{shafieezadeh2023constrained}. The above} observation motivates the need for a comprehensive study of the mixed-integer set
\begin{equation}\label{eq:def-Q}
\hspace{-0.9em}	\set Q\defeq\left\{ (t,x,z)\in\R^{n+1}\times\{0,1\}^n:\;\begin{aligned}
		&t\ge f(x), x_i\ge 0\;\forall i\in\set I_+,\\
		& z\in \{0,1\}^n,x_i(1-z_i)=0\;\forall i\in[n]
	\end{aligned}\right\},
\end{equation}
where $f(x)=g(Ax)+{c}^\top{x}$, $g:\R^k\to\R$ is a proper closed convex function, $A$ is a $k\times n$ matrix, $c\in\R^n$, and $\set I_+\subseteq [n]$ is the subset of variables restricted to be non-negative.

Disjunctive programming is a powerful modeling tool to represent a nonconvex optimization problem in the form of disjunctions, especially when binary variables are introduced to encode the logical conditions such as sudden changes, either/or decisions, implications, etc. \cite{balas2018disjunctive}. The theory of linear disjunctive programming {was} first pioneered by Egon Balas in 1970s \cite{balas1998disjunctive, balas1979disjunctive, balas1985disjunctive, balas1989sequential}, and later extended to the nonlinear case \cite{Ceria1999, stubbs1999branch, lee2000new, grossmann2002review,bernal2021convex,bonami2011lift,kilinc2010effective,kilincc2015two,modaresi2016intersection,ccezik2005cuts,yildiz2016low}. Once a mixed-integer set is modeled as a collection of disjunctive sets, its convex hull can be described easily as a Minkowski sum of the scaled convex sets with each obtained by creating a copy of original variables. Such extended formulations are the strongest possible convex relaxations of the mixed-integer set under study. However, a potential downside of such formulations is that, often, the number of additional variables required in the description of the convex hull is exponential in the number of binary variables. Thus, a direct application of disjunctive programming as mentioned can be unfavorable in practice.

A possible approach to implement the extended formulation induced by disjunctive programming in an efficient way is to reduce the number of additional variables introduced in the model without diminishing the relaxation strength. In principle, this goal can be achieved by means of Fourier-Motzkin elimination \cite{dantzig1972fourier}. This method is practical if the set under study can be naturally expressed using few disjunctions, e.g., to describe piecewise linear functions \cite{anderson2019strong,han2021single}, {involves} few binary variables \cite{atamturk2018strong,frangioni2019decompositions}, or {is} separable \cite{gunluk2010perspective}. However, projecting out variables can be very challenging even if $n=2$ \cite{anstreicher2021quadratic,de2022explicit}, not to mention in a high-dimensional setting.
 Regarding the nonlinear set with indicator variables $\set Q$, in the simplest case where the MINLP is separable, it is known that the convex hull can be described by the perspective function of the objective function supplemented by the constraints defining the feasible region \cite{akturk2009strong,dong2015regularization,dong2013,frangioni2006perspective, gunluk2010perspective,
	hijazi2012mixed,wu2017quadratic, bertsimas2020mixed}. 

Two papers \cite{wei2020ideal} and \cite{atamturk2020supermodularity} are closely related to our work. The mixed-integer sets studied in both can be viewed as a special realization of $\set Q$ in \eqref{eq:def-Q}. \citet{wei2020ideal} characterize the closure of the convex hull if $\set Q$ when the defining function $g$ is a univariate function and the continuous variables are free, i.e. $\set I_+=\emptyset$. \citet{atamturk2020supermodularity} study the setting with sign-constrained continuous variables, and univariate quadratic functions $g$: in such cases, the sign-restrictions resulted in more involved structures for the closure of the convex hull of $\set Q$.
\subsection*{Contributions and outline}
In this paper we show how to construct extended formulations of set $\set Q$, requiring at most $\mathcal{O}(n^k)$ copies of the variables. In particular, if the {dimension} $k$ is small, then the resulting formulations are indeed much smaller than those resulting from a natural application of disjunctive programming. Moreover, for the special case of $k=1$, we show how we are able to recover and improve existing results in the literature, either by providing smaller extended formulations (linear in $n$), or by providing convexifications in the original space of variables for a more general class of functions {(arbitrarily convex and not necessarily differentiable)}. {When applied to a class of sparse signal denoising and outlier detection problems, the new formulations proposed in this paper significantly outperforms the natural one whose relaxation yields a trivial lower bound of $0$. It also turns out that as $k$ increases, the resulting formulation becomes even more powerful.}
 
The rest of the paper is organized as follows. In Section~\ref{sec:cvx-analysis}, we provide relevant background and the main result of the paper: a compact extended formulation of $\set Q$. In Section~\ref{sec:rank-one}, using the results in Section~\ref{sec:cvx-analysis}, we derive the explicit form of the convex hull $\set Q$ {for $k=1$} in the original {space} of variables. 
 In Section~\ref{sec:complexity}, we present the complexity results showing that tractable convexifications of $\set{Q}$ are unlikely if additional constraints are imposed on the continuous variables. {In Section~\ref{sec:experiments}, we present the numerical experiments.}
 Finally, in Section~\ref{sec:conclusions}, we conclude the paper.

\section{A convex analysis perspective on convexification}\label{sec:cvx-analysis}
In this section, we first introduce necessary preliminaries in convex analysis and notations adopted in this paper. After that we present our main results and their connections with previous works in literature.
\subsection{Notations and preliminaries} Throughout this paper, we assume $f(x)$ {is} a proper closed convex function from $\R^n$ to $\R\cup\{+\infty\}$. If $f(x)<+\infty$ for all $x$, $f(x)$ is called \emph{finite}. We denote the \emph{effective domain} of $f$ by $\dom(f)$ and  the convex \emph{conjugate} of $f(\cdot)$ by $f^*(\cdot)$ which is defined as $$f^*(\alpha)\defeq\max_{x}\;{\alpha}^\top{x}-f(x).$$ The \emph{perspective function} of $f(\cdot)$ is defined as \[f^\pi(x,\lambda)\defeq\begin{cases}
	\lambda f\left(\frac{x}{\lambda}\right)&\text{if } \lambda\ge0\\
	+\infty&\text{o.w.},
\end{cases}\]
where $0f(x/0)$ is interpreted as the \emph{recession function} of $f$ at $x$. By this definition, {the} function $f^\pi(\cdot)$ is homogeneous, closed and convex; see Section~8 in \cite{rockafellar1970convex}. Moreover, we borrow the term \emph{rank} which is normally defined for an affine mapping and extend its definition to a general nonlinear convex function.
\begin{definition}[Rank of convex functions]
	Given a proper closed convex function $f(x)$, the rank\footnote{The definition of $\rank(f)$ is different from the one adopted in classical convex analysis (see Section~8 in \cite{rockafellar1970convex}). If $\dom(f)$ is full-dimensional, the two definitions coincide.} {of $f$, denoted $\rank(f)$,} is defined as the smallest integer $k$ such that $f(x)$ can be expressed in the form $f(x)=g(Ax)+{c}^\top{x}$ for some closed convex function $g(\cdot)$, a $k\times n$ matrix $A$ and a vector $c\in\R^n$.
\end{definition}
For example, the rank of an affine function $f(x)={c}^\top{x}$ is simply $0$. The rank of a convex quadratic function $f(x)=x^\top A{x}$ coincides with $\rank(A)$, where $A\succeq 0$ is a positive semidefinite (PSD) matrix. 

We let $\ones$ be the vector of all ones (whose dimension can be inferred from the context). {For any vector $y\in \R^n$ and index set $\set J\subseteq [n]$, we denote $y(\set J)=\sum_{j\in\set J}y_j$.} For a set $\set S$, we denote the convex hull of $S$ as $\text{conv}(\set S)$, and its closure as $\conv(\set S)$. For any scalar ${\lambda}\ge0$ and two generic sets $\set S_1$, $\set S_2$ in a proper Euclidean space, we define $\lambda \set S_1\defeq\{\lambda x:x\in\set S_1\}$ and $\set S_1+\set S_2\defeq\{x^1+x^2:x^1\in\set S_1,x^2\in\set S_2 \}$ is the Minkowski sum. We denote the \emph{indicator function} of $\set S$ by $\delta(x;\set S)$, which is defined as $\delta(x;\set S)=0$ if $x\in \set S$ and $\delta(x;\set S)=+\infty$ otherwise. By the above notations, the convex conjugate of $\delta(x;\set S)$ is $$\delta^*(\alpha;\set S)=\max_x\; {\alpha}^\top{x}-\delta(x;\set S)= \max_{x\in\set S}\;{\alpha}^\top{x},$$ which is known as the \emph{support function} of $\set S$. To be consistent with the definition of $0f(x/0)$, we interpret $0\set S$ as the \emph{recession cone} of $\set S$ throughout this paper. The derivation of the succeeding work relies on the following one-to-one correspondence between closed convex sets and support functions, whose proof can be found in classical books of convex analysis, e.g. Section~13 in \cite{rockafellar1970convex} and Chapter~C.2 in \cite{hiriart2004fundamentals}.
\begin{proposition}\label{prop:outer-characterization}
	Given a set $\set S$ and a closed convex set $\set T$, $\set T=\conv(\set S)$ if and only if $\delta^*(\cdot;\set T)=\delta^*(\cdot;\set S)$.
\end{proposition}

For convenience, we repeat the set of interest:
\[ \set Q\defeq\left\{ (t,x,z)\in\R^{n+1}\times\{0,1\}^n:\quad\begin{aligned}
	&t\ge f(x), x_i\ge 0\;\forall i\in\set I_+,\\
	&x_i(1-z_i)=0\;\forall i\in[n]
\end{aligned}\right\}, \]
where $f(x)=g(Ax)+{c}^\top{x}$ is a proper closed convex function.

 Note that if the complementary constraints $x_i(1-z_i)=0$ are removed from $\set Q$, then $x$ and $z$ are decoupled and $\conv(\set Q)$ reduces simply to $\set X\times[0,1]^n$, where $$\set X\defeq\left\{ (t,x)\in\R^{n+1}:t\ge f(x), x_i\ge 0\;\forall i\in\set I_+ \right\}.$$   For the purpose of decomposing $\conv(\set Q)$, for any $\set I\subseteq[n]$, define the following sets: 
 \begin{equation}\label{eq:def-VI}
 	\begin{aligned}
 		&\set X(\set I)\defeq\set X\cap\{(t,x):x_i=0\;\forall i\notin \set I \},\\
 		&\set Z(\set I)\defeq\left\{z\in\{0,1\}^n:z_i=1\;\forall i\in \set I \right\},\\
 		&\set V(\set I)\defeq\set X(\set I)\times\set Z(\set I). 
 	\end{aligned}
 \end{equation}
Notice that for any $(t,x,z)\in\set V(\set I)$, either $x_i=0$ or $z_i=1\;\forall i\in[n]$. Therefore, $x_i(1-z_i)=0\;\forall i\in[n]$ and thus $\set V(\set I)\subseteq \set Q$. Furthermore, $\set Q$ can be expressed as a disjunction
\begin{equation}\label{eq:naive-disjunction}
	\set Q=\bigcup_{\set I\subseteq[n]} \set V(\set I).
\end{equation}
{Standard} disjunctive programming techniques {for convexification of $\mathcal{Q}$} are based on \eqref{eq:naive-disjunction}, creating copies of variables for every $ I\subseteq[n]$, resulting in an exponential number of variables.

Finally, define  
\begin{equation}\label{eq:def-R}
	\set R\defeq \{(t,x,z):t\ge0,\;Ax=0, x_i\ge0,\;\forall i\in\set I_+,\;z_i=1,\;\forall i\in[n] \},
\end{equation}
which is a  closed convex set.
 For any $y\in\R^n$, we denote the \emph{support} of $y$ by $\support(y)\defeq\{i:y_i\neq0  \}$

\subsection{Convex hull characterization}
In this section, we {aim to} characterize {$\conv(\set Q)$}. We first show that if $f(x)$ is homogeneous, under mild conditions, $\conv(\set Q)$ is simply the natural relaxation of $\set Q$. 
\begin{proposition}\label{prop:homogeneous}
	If $f(x)$ is a homogeneous function, then $\conv(\set Q)=\set X\times[0,1]^n$.
\end{proposition} 
\begin{proof}
	Since $f$ is homogeneous, each set $\set X(\set I)$ is a closed convex cone. Moreover, $\set I_1\subseteq \set I_2$ implies that $\set X(\set I_1)\subseteq\set X(\set I_2)$, which further implies that 
	\begin{equation}\label{eq:sum-XI}
	\set X(\set I_2)\subseteq	\set X(\set I_1)+\set X(\set I_2)\subseteq \set X(\set I_2)+\set X(\set I_2)=\set X(\set I_2),
	\end{equation}
where the last equality results from that $\set X(\set I_2)$ is conic. Hence, equality holds throughout \eqref{eq:sum-XI}. By \eqref{eq:naive-disjunction},
	\begin{align*}
		&\conv(\set Q)=\conv\left( \bigcup_{\set I\subseteq[n]} \set V(\set I) \right)\\
		=&\conv\left(\; \smashoperator{\bigcup_{\set I\subseteq[n]}} \set X(\set I)\times\text{conv}(\set Z(\set I)) \right)\tag{\small $\set X(\set I)$ is convex; $\set{X},\set{Z}$ are decoupled}\\
		=&\bigcup_{\lambda\in\R_+^{2^n}:\ones^\top \lambda=1}\;\sum_{\set I\subseteq[n]}\lambda_{\set I}\left(\set X(\set I)\times\text{conv}(\set Z(\set I))\right)\\
		=&\bigcup_{\lambda\in\R_+^{2^{n}}:\ones^\top \lambda=1}\;\set X([n])\times\sum_{\set I\subseteq[n]}\lambda_{\set I} \text{conv}(\set Z(\set I))\tag{$\set X(\set I)$ is conic and \eqref{eq:sum-XI}}\\
		=&\set X\times\bigcup_{\lambda\in\R_+^{2^n}:\ones^\top \lambda=1}\;\sum_{\set I\subseteq[n]}\lambda_{\set I} \text{conv}(\set Z(\set I))\tag{$\set X$=$\set X([n])$}\\
		=&\set X\times [0,1]^n.
	\end{align*}
\end{proof}
Proposition~\ref{prop:homogeneous} generalizes Proposition~1 of \cite{gomez2021strong}{, from the fact that $f$ is the $\ell_2$ norm}. Next, we present the main result of the paper, characterizing $\conv(\set Q)$ without the assumption of homogeneity. In particular, we show that $\conv(\set Q)$ can be constructed from substantially {fewer} disjunctions than those given in~\eqref{eq:naive-disjunction}.
\begin{theorem}\label{thm:main}
	Assume $\rank(f)\le k$ and $f(0)=0$. Then
	 \[ \conv(\set Q)=\conv\left( \left(\bigcup_{\set I:|\set I|\le k}\set V({\set I})\right)\cup \set R \right){,}\]
	 {where $\set V(\set I)$ and $\set R$ are defined in \eqref{eq:def-VI} and \eqref{eq:def-R} respectively.}
	  Moreover, if $\left\{x\in\R^n:Ax=0,\; x_i\ge0\;\forall i\in\set I_+ \right\}=\{0\}$, then $\set R$ can be removed from the disjunction.
\end{theorem}
Informally, $\bigcup_{\set I:|\set I|\le k}\set V(I)$ and $\set R$ correspond to the ``extreme points" and ``extreme rays" of $\conv(\set Q)$, respectively. From \eqref{eq:naive-disjunction}, $\set Q$ is a disjunction of exponentially many pieces of $\set V(\set I)$. However, given a low-rank function $f$, Theorem~\ref{thm:main} states that $\conv(\set Q)$ can be generated (using disjunctive programming) from a much smaller number of sets $\set V(\set I)$.  We also remark that condition $f(0)=0$ plays a minor role in the derivation of Theorem~\ref{thm:main}. If $0\in\dom(f)$ but $f(0)\neq 0$, one can study $f(x)-f(0)$. 
\begin{proof}[Proof of Theorem~\ref{thm:main}]
	Denote $\set S=\smashoperator{\bigcup\limits_{\set I:|\set I|\le k}}\set V(I)\cup \set R$. Since $\set V(\set I)\subseteq \set Q$ for all $\set I$, and $\set{R}\subseteq \set Q$, we find that $\set S\subseteq \set Q$.  Thus, $\delta^*(\cdot;\conv(\set S))=\delta^*(\cdot;\set S)\le\delta^*(\cdot;\set Q). $	Due to Proposition~\ref{prop:outer-characterization}, it remains to prove the opposite direction, namely that $\delta^*(\cdot;\set S)\ge\delta^*(\cdot;\set Q)$. 
	
	Since $\rank(f)\le k$, there exists $g(\cdot)$, $A\in\R^{k\times n}$ and $c$ such that $f(x)=g(Ax)+{c}^\top{x}$. Taking any $\alpha=(\alpha_t,\alpha_x,\alpha_z)\in \R\times \R^n\times \R^n$, if $\alpha_t>0$, $\delta^*(\cdot;\set S)=\delta^*(\cdot;\set Q)=+\infty$ because $t$ is unbounded from above. We now assume $\alpha_t\le 0$, and define $$\ell(x,z)\defeq\alpha_tg(Ax)+\left(\alpha_x+\alpha_t c\right)^\top{x}+\alpha_z^\top{z}.$$ Then $\delta^*(\alpha;\set S)=\max\{\ell(x,z):(t,x,z)\in\set S\}$ and $\delta^*(\alpha;\set Q)=\max\{\ell(x,z):(t,x,z)\in\set Q\}$. Given any fixed $(\bar t,\bar x,\bar z)\in\set Q$, consider the linear program
	\begin{equation}\label{eq:main-lp}
		\begin{aligned}
			v^*\defeq\max_{x\in\R^n}\;&\left(\alpha_x+\alpha_tc\right)^\top{x}\\
			\text{s.t. }& Ax=A\bar x\\
			&\bar x_ix_i\ge 0&\forall i\in\support(\bar x)\subseteq \support(\bar z)\\
			&x_i=0&\forall i\in [n]\backslash\support(\bar x).
		\end{aligned}
	\end{equation}
Note that linear program \eqref{eq:main-lp} is always feasible setting $x=\bar x$. Moreover, every feasible solution (or direction) $x$ of \eqref{eq:main-lp} satisfies $x_i(1-\bar z_i)=0$.

 If $v^*=+\infty$, then there exists a feasible direction $d$ such that $Ad=0$ and $\left(\alpha_x+\alpha_tc\right)^\top{d}>0$. It implies that for any $r\ge0$, $f(rd)=g(rAd)+rc^\top{d}=g(0)+r{c}^\top{d}=r{c}^\top{d}$. Hence, $(r{c}^\top{d},rd,\ones)\in\set R\subseteq Q$. Furthermore, 
\begin{align*}
	\ell(rd,\ones) =&\left(\alpha_x+\alpha_tc\right)^\top{d}r+{\alpha^\top_z}{\ones}
	\rightarrow+\infty,\tag{as $r\rightarrow+\infty$}
\end{align*}  which implies $\delta^*(\alpha;\set S)=\delta^*(\alpha;\set Q)=+\infty$.

 If $v^*$ is finite, there exists an optimal solution $x^*$ to \eqref{eq:main-lp}. Moreover, $x^*$ can be taken as an extreme point of the feasible region of \eqref{eq:main-lp} which is a pointed polytope. It implies that $n$ linearly independent constraints must be active at $x^*$. Since $\rank(A)\le k$, at least $n-k$ constraints of the form $\bar x_i x_i\geq 0$ {and $x_i=0$ in \eqref{eq:main-lp}} hold at equality. Namely, $x^*$  satisfies $|\support(x^*)|\le k$. Since $A\bar x=A x^*$, we can define
 \begin{align*}
 	&t^*\defeq\bar t+{c^\top}{{(x^*-\bar x)}}\\
 	\ge& f(\bar x)+{c^\top}{{(x^*-\bar x)}}\\
 	=&g(A\bar x)+{c^\top}{\bar x}+{c^\top}{{(x^*-\bar x)}}\\
 	=&g(Ax^*)+{c^\top}{x^*}=f(x^*)
 \end{align*} Setting $\set I=\support(x^*)$, one can deduce that $( t^*, x^*)\in\set X(\set I)$, and $\bar z\in\set Z(\set I)$ because  $\support(x^*)\subseteq \support(\bar z)$.   Thus, $(t^*, x^*,\bar z)\in \set V(\set I)\subseteq \set S\subseteq \set Q$, and $\left(\alpha_x+\alpha_tc\right)^\top{x^*}\ge \left(\alpha_x+\alpha_tc\right)^\top{\bar x}$ {implies} that  $\ell(x^*, \bar z)\ge \ell(\bar x,\bar z)$. Namely, for an arbitrary point $(\bar t,\bar x,\bar z)\in \set Q$, there always exists a point in $\set S$ with a superior {or equal} objective value of $\ell(\cdot)$. Therefore, $\delta^*(\alpha;\set S)\ge\delta^*(\alpha;\set Q)$, completing the proof of the main conclusion. 
 
 The last statement of the theorem follows since if 
$\{x\in\R^n:Ax=0,\; x_i\ge0\;\forall i\in\set I_+ \}=\{0\}$, then the feasible region of \eqref{eq:main-lp} is bounded and thus $v^*$ is always finite.
\end{proof}

Using the disjunctive representation of Theorem~\ref{thm:main} and usual disjunctive programming techniques \cite{Ceria1999}, one can immediately obtain extended formulations requiring at most $\mathcal{O}(n^k)$ copies of the variables. Moreover, it is often easy to project out some of the additional variables, resulting in formulations with significantly {fewer} variables or, in some cases, formulations in the original space of variables. We illustrate these concepts in the next section with $k=1$ {and in Section~\ref{sec:experiments} for $k=2$}.

\section{Rank-one convexification}\label{sec:rank-one}
In this section, we show how to use Theorem~\ref{thm:main} and disjunctive programming to derive convexifications for rank-one functions. In particular, throughout this section, we make the following assumptions:

\begin{assumption}\label{assum:r1}
	Function $f$ is given by $f{(x)}=g\left(\sum_{i=1}^na_ix_i\right)+{c}^\top{x}$, where $g$ is a finite one-dimensional function, $f(0)=g(0)=0$, and $a_i\neq0\; \forall i\in n$. For simplicity, we also assume $c=0$.
\end{assumption}
Note that the value specification of $g(0)$ and $c$ is assumed without loss of generality in Assumption~\ref{assum:r1}. Indeed, for a general finite rank-one convex function $f$, $(t,x,z)\in\conv(\set Q)$ if and only if $(t-c^\top x-f(0), x, z)\in\conv(\set Q_{\tilde f})$, where $\tilde f(x)\defeq f(x)-c^\top x-f(0)$ and $\set Q_{\tilde f}$ is the corresponding mixed-integer set associated with the variant $\tilde f$. The goal of this section is to characterize the special case of $\set Q$ given by
\[ \setRankOne\defeq\left\{ (t,x,z)\in\R^{n+1}\times\{0,1\}^n:\quad\begin{aligned}
	&t\ge f(x)=g\left(\sum_{i=1}^na_ix_i\right), \\
	&x_i\ge 0\;\forall i\in\set I_+,\;x_i(1-z_i)=0\;\forall i\in[n]
\end{aligned}\right\}. \] 

First in Section~\ref{sec:extended} we derive an extended formulation {$\conv(\setRankOne)$} with a {linear} number of additional variables. Then in Section~\ref{sec:original_free} we project out the additional variables for the case with free continuous variables, and recover the results of \cite{wei2020ideal}. Similarly, in Section~\ref{sec:original_nonneg}, we provide the description of the convex hull of cases with non-negative continuous variables in the original space of variables, generalizing the results of \cite{atamturk2020supermodularity} and \cite{atamturk2018strong} to general (not necessarily quadratic) functions $g$. We also show that the extended formulation proposed in this paper is more amenable to implementation that the inequalities proposed in \cite{atamturk2020supermodularity}. 
\subsection{Extended formulation of $\conv({\setRankOne})$}\label{sec:extended}
We first discuss the compact extended formulation of $\conv({\setRankOne})$ that can be obtained directly for the disjunctive representation given in Theorem~\ref{thm:main}, using $2n$ additional variables. 
\begin{proposition}\label{prop:rank-one} Under Assumption~\ref{assum:r1}, 
	$(t,x,z)\in\conv({\setRankOne})$ if and only if there exists $\lambda,\tau\in \R^n$ such that the inequality system 
	\begin{equation}\label{eq:rank-one}
		\begin{aligned}
			&t\ge \sum_{i=1}^ng^\pi(a_i(x_i-\tau_i),\lambda_i),\\
			&a^\top \tau=0,\;0\le \tau_i\le x_i\;\forall i\in\set I_+,\\
			&\lambda_i\le z_i\le 1\;\forall i\in[n],\\
			&\lambda\ge 0,\;\sum_{i=1}^n \lambda_i\le1
		\end{aligned}
	\end{equation}
	is satisfied.
\end{proposition}
\begin{proof}
 By Theorem~\ref{thm:main}, 
	\begin{align*}
		&\conv({\setRankOne})=\conv\left( \set V(\emptyset)\bigcup_{i\in[n]}\set V(\{i\})\cup\set R \right)\\
		=&\bigcup_{\lambda\in\R_+^{n+2}:\ones^\top \lambda=1}\left(\lambda_0\conv(\set V(\emptyset))+\sum_{i=1}^n\lambda_i\set \conv({\set V(\{i\})})+\lambda_{[n]}\set R\right),\\
	\end{align*}
where 
\begin{align*}
	&\lambda_0\conv(\set V(\emptyset))=\{ (t_0,x^0,z^0):t_0\ge0,\;x^0=0,\;0\le z^0\le \lambda_0 \},\\
	&\lambda_i\conv(\set V(\{i\}))=\left\{(t_i,x^i,z^i):\begin{aligned}
		&t_i\ge g^\pi(a_ix^i_i,\lambda_i),\;z^i_i=\lambda_i,\;0\le z^i_j\le\lambda_i\;\forall j\neq i,\\
		&x^i_i\ge0 \text{ if } i\in\set I_+,\;x_j^i=0\;\forall j\neq i
	\end{aligned} \right\},\\
&\lambda_{[n]}\set R=\left\{\left(t_{[n]},\tau,z^{[n]}\right):t_{[n]}\ge0,\;a^\top \tau=0,\;\tau_i\ge0\;\forall i\in\set I_+,\;z^{[n]}=\lambda_{[n]}\ones\right\}.
\end{align*}
It follows that $(t,x,z)\in\conv({\setRankOne})$ if and only if the following inequality system has a solution:
\begin{subequations}
	\begin{align*}
		&t\ge \sum_{i=1}^ng^\pi(a_ix^i_i,\lambda_i),\;a^\top \tau=0,\\
		&x_i=x_i^i+\tau_i&\forall i\in[n],\numberthis\label{eq:reduce-xii}\\
		&x_i^i\ge0,\;\tau_i\ge0&\forall i\in\set I_+,\\
		&z_j= \sum_{i=1}^nz_j^i+\lambda_{[n]}+z_j^0&\forall j\in[n],\\
		&z_i^i=\lambda_i,\;0\le z_i^0\le \lambda_0&\forall i\in[n],\numberthis\label{eq:reduce-zii}\\
		&0\le z_j^i\le \lambda_i&\forall j\neq i\in[n],\\
		&\lambda\ge 0,\;\lambda_0+\sum_{i=1}^n \lambda_i+\lambda_{[n]}=1.
	\end{align*}
\end{subequations}
We now show how to {simplify} the above inequality system step by step. First, we can substitute out $x^i_i$ and $z^i_i$ using \eqref{eq:reduce-xii} and \eqref{eq:reduce-zii}, obtaining the system
\begin{subequations}
	\begin{align*}
		&t\ge \sum_{i=1}^ng^\pi(a_i(x_i-\tau_i),\lambda_i),\;a^\top \tau=0,\\\
		&0\le \tau_i\le x_i&\forall i\in\set I_+,\\
		&z_j-\lambda_j-\lambda_{[n]}=\sum_{i\in[n]:i\neq j}z^i_j+z_j^0&\forall j\in[n],\numberthis\label{eq:reduce-zij-1}\\
		&0\le z_j^i\le \lambda_i&\forall [n]\ni i\neq  j\in[n]\cup\{0\},\numberthis\label{eq:reduce-zij-2}\\
		&\lambda\ge 0,\;\lambda_0+\sum_{i=1}^n \lambda_i+\lambda_{[n]}=1.\numberthis\label{eq:reduce-lambda0}
	\end{align*}
\end{subequations}
Next, we can substitute out $\sum_{i\in[n]:i\neq j}z^i_j+z_j^0$ in \eqref{eq:reduce-zij-1} using the bounds \eqref{eq:reduce-zij-2}. Doing so, \eqref{eq:reduce-zij-1} reduces to
\begin{align*}
	&0\le z_j-\lambda_j-\lambda_{[n]}\le \sum_{i\in[n]:i\neq j}\lambda_i+\lambda_0=1-\lambda_j-\lambda_{[n]},\\
	\Leftrightarrow& \lambda_{[n]}\le z_j-\lambda_j,\;\; z_j\le 1,
\end{align*}
where the equality results from \eqref{eq:reduce-lambda0}. We deduce that the system of inequalities reduces to
\begin{subequations}
	\begin{align*}
		&t\ge \sum_{i=1}^ng^\pi(a_i(x_i-\tau_i),\lambda_i),\; a^\top \tau=0,\\\
		&0\le \tau_i\le x_i &\forall i\in\set I_+,\\
		&z_j\le 1, \;\lambda_{[n]}\le z_j-\lambda_j&\forall j\in[n],\\
		&\lambda\ge 0,\;\lambda_0+\sum_{i=1}^n \lambda_i+\lambda_{[n]}=1.
	\end{align*}
\end{subequations}
Formulation~\eqref{eq:rank-one} follows from using Fourier-Motzkin elimination to project out $\lambda_0$ and $\lambda_{[n]}$, replacing them with 0 and changing the last equality to an inequality.
\end{proof}

 In addition, if $\set I_+=[n]$ and $a_i>0\;\forall i\in[n]$, then $a^\top\tau=0$ and $\tau\ge0$ imply $\tau=0$ in \eqref{eq:rank-one}. Therefore, we deduce the following corollary for this special case.
\begin{corollary}\label{cor:rank-one-pos}
	Under Assumption~\ref{assum:r1}, $\set I_+=[n]$ and $a_i>0\;\forall i\in[n]$, $(t,x,z)\in\conv({\setRankOne})$ if and only if there exists $\lambda\in \R^n$  such that the inequality system 
	\begin{equation}
		\begin{aligned}
			&t\ge \sum_{i=1}^ng^\pi(a_ix_i,\lambda_i),\\
			&\lambda_i\le z_i\le 1\;\forall i\in[n],\\
			&\lambda\ge 0,\;\sum_{i=1}^n \lambda_i\le1
		\end{aligned}
	\end{equation}
is satisfied.
\end{corollary}

Naturally, the extended formulations obtained from Proposition~\ref{prop:rank-one} and Corollary~\ref{cor:rank-one-pos}, requiring $\mathcal{O}(n)$ additional variables, are substantially more compact than extended formulations obtained from the disjunction \eqref{eq:naive-disjunction}, which require $\mathcal{O}(n^{2^n})$ variables. 
\subsection{Explicit form of $\conv({\setRankOne})$ with unconstrained continuous variables}\label{sec:original_free}
When $x$ is unconstrained, i.e., $\set I_+=\emptyset$, the explicit form of $\conv({\setRankOne})$ in the original space of variables was first established in \cite{atamturk2019rank} for quadratic functions, and later generalized in \cite{wei2020ideal} to general rank-one functions. In Proposition~\ref{prop:ran-one-explicit} below, we present a short proof on how to recover the aforementioned results, starting from Proposition~\ref{prop:rank-one}. First, we need the following property on the monotonicity of the perspective function.
\begin{lemma}\label{lem:monotonicity-persp}
		Assume $g(v)$ is a convex function over $\R$ with $g(0)=0$. Then $g^\pi(v,\lambda)$ is a nonincreasing function on $\lambda>0$ for fixed $v$.
\end{lemma}
\begin{proof}
	Since $g$ is convex, for $\forall v_1,v_2$,
	 $\frac{g(v_1)-g(v_2)}{v_1-v_2}$ is nondecreasing with respect to $v_1$. Taking $v_1=v/\lambda$ and $v_2=0$, since $g(0)=0$, one can deduce that
	 \[g^\pi(v,\lambda)=\lambda g\left(\frac{v}{\lambda}\right)= v\frac{g\left(v\frac{1}{\lambda}-0\right)}{v\frac{1}{\lambda}-0} \]
	 is nondecreasing with respect to $1/\lambda$, i.e. nonincreasing with respect to $\lambda$.
\end{proof}
\begin{proposition}[\citet{wei2020ideal}]\label{prop:ran-one-explicit}
Under Assumption~\ref{assum:r1}, if additionally $\set I_+=\emptyset$, then
\[\conv({\setRankOne})=\left\{ (t,x,z):t\ge g^\pi\left( a^\top x,\,\min\left\{1,\,\ones^\top z \right\} \right),\;0\le z\le \ones \right\}. \]
\end{proposition}
\begin{proof}
	Without loss of generality, we assume $a=\ones$; otherwise, we can scale $x_i$ by $a_i$. We first eliminate $\tau$ from \eqref{eq:rank-one}. From the convexity of $g(\cdot)$, we find that  
	\begin{align*}
		&\sum_{i\in[n]}g^\pi(x_i-\tau_i,\lambda_i)=\sum_{i\in[n]}\lambda_ig\left(\frac{x_i-\tau_i}{\lambda_i}\right)\\
		=&\left(\sum_{j\in[n]}\lambda_j\right)\left(\sum_{i\in[n]}\frac{\lambda_i}{\sum_{j\in[n]}\lambda_j}g\left( \frac{x_i-\tau_i}{\lambda_i} \right) \right)\\
		\ge&\left(\sum_{j\in[n]}\lambda_j\right)g\left(\frac{\sum_{i\in[n]}x_i-\sum_{i\in [n]}\tau_i}{\sum_{j\in[n]}\lambda_j} \right)\\
		=&\left(\sum_{j\in[n]}\lambda_j\right)g\left(\frac{\sum_{i\in[n]}x_i}{\sum_{j\in[n]}\lambda_j} \right)=g^\pi\left(\sum_{i\in[n]}x_i,\sum_{i\in[n]}\lambda_i\right),\tag{$\ones^\top \tau=0$}
	\end{align*}
where the inequality holds at equality if there exists some common ratio $r$ such that $x_i-\tau_i=\lambda_i r$. Moreover, this ratio does exist by setting $r=\ones^\top x/\ones^\top \lambda$ and $\tau=x-r\lambda$ for all $i\in [n]$ --it can be verified directly that $\ones^\top \tau=0$. Thus, the above lower bound can be attained for all $x$. Hence, \eqref{eq:rank-one} reduces to
\begin{subequations}
	\begin{align}
		&t\ge g^\pi(\ones^\top x,\ones^\top \lambda)\nonumber\\
		&0\le \lambda_i\le z_i\le 1\;\forall i\in[n]\label{eq:reduce-lambda-1}\\
		&\sum_{i\in[n]}\lambda_i\le 1\label{eq:reduce-lambda-2}.
	\end{align}
\end{subequations}
Since for fixed $v$, $g^\pi(v,s)$ is  non-increasing with respect to $s\in\R_+$, projecting out $\lambda$ amounts to computing the maximum of $\ones^\top \lambda$, that is, solving the linear program $\max_\lambda\;\{\ones^\top \lambda:\eqref{eq:reduce-lambda-1}\text{ and }\eqref{eq:reduce-lambda-2} \}.$ Summing up \eqref{eq:reduce-lambda-1} over all $i\in[n]$ and combining it with \eqref{eq:reduce-lambda-2}, we deduce that $\ones^\top \lambda\le\min\{\ones^\top z,\;1  \}$. It remains to show this upper bound is tight. If $\sum_{i\in[n]}z_i\le 1$, one can set $\lambda_i=z_i\;\forall i\in[n]$. Now assume $\sum_{i\in[n]}z_i> 1$. Let $m$ be the index such that $\sum_{i=1}^{m-1}z_i\le 1$ and $\sum_{i=1}^mz_i>1$. Set 
\begin{align*}
	\lambda_i=\begin{cases}
		z_i&\text{if }i<m\\
		1-\sum_{i=1}^{m}z_i&\text{if }i=m\\
		0&\text{if } i>m.
	\end{cases}
\end{align*} It can be verified directly that this solution is feasible and $\ones^\top\lambda=1$. The conclusion follows.
\end{proof}
\subsection{Explicit form of $\conv({\setRankOne})$ with nonnegative continuous variables}\label{sec:original_nonneg}
In this section, we aim to derive the explicit form of $\conv({\setRankOne})$ when $\set I_+=[n]$ {in the original space of variables. In other words, we show how to project out the additional variables $\lambda$ and $\tau$ from formulation \eqref{eq:rank-one}}. A description of this set is known for the quadratic case \cite{atamturk2020supermodularity} only. We now derive it for the general case. When specialized to bivariate rank-one functions, it also generalizes the main result of \cite{atamturk2018strong} where the research objective is a non-separable bivariate quadratic.   

 {Slightly abusing the notation $a$, we observe} that function $f$ can be written {in the form}
$$f(x)=g\left(\sum_{i\in\set N_+}a_ix_i-\sum_{i\in\set N_-}a_ix_i\right),$$
where $a_i>0\;\forall i\in[n]$, and  $\set N_+\cup\set N_-$ is a partition of $[n]$.
{Theorem~\ref{thm:explicit-pos} below contains the main result of this section. It gives necessary and sufficient conditions under which an arbitrary point $(t,x,z)\in \R^{2n+1}$ belongs the closure of the convex hull of $\setRankOne$. Throughout,}
the $\min/\max$ over an empty set is taken to be $+\infty/{-}\infty$, respectively.

 \begin{theorem}\label{thm:explicit-pos}
 	Under Assumption~\ref{assum:r1} and $\set I_+=[n]$, for all $(t,x,z){\in \R^{2n+1}}$ such that $x\ge 0, 0\le z\le \ones$, the following statements hold:
	\begin{itemize}
		\item If $\sum_{i\in \set N_+}a_ix_i>\sum_{i\in\set N_-}a_ix_i$ and there exists a partition $\set L\cup \set M\cup\set U$ of $\set N_+\cap\support(x)\cap\support(z)$ such that 
		\begin{align*}
			&1-\smashoperator{\sum_{i\in\set M\cup\set U}}z_i> 0,\quad\max_{i\in\set L}\frac{a_ix_i}{z_i}<\frac{\sum_{i\in \set L}a_ix_i}{1\!\!-\!\!\sum_{i\in\set M\cup\set U}z_i}\le\min_{i\in\set M}\frac{a_ix_i}{z_i}\\
			&\max_{i\in\set M}\frac{a_ix_i}{z_i}\le \frac{\sum_{i\in\set N_+}a_ix_i-\sum_{i\in\set L\cup\set M\cup\set N_-}a_ix_i}{\sum_{i\in\set U}z_i}<\min_{i\in\set U}\frac{a_ix_i}{z_i},
		\end{align*}
	then $(t,x,z)\in\conv({\setRankOne})$ if and only if 
	\begin{align}\label{eq:convex-envolope}
		\hspace{\leftmargin}
		\begin{multlined}[b][.8\textwidth]
			t\ge g^\pi\left( {\sum_{i\in \set L}a_ix_i},{1\!-\!\smashoperator{\sum_{i\in\set M\cup\set U}}z_i} \right)+\sum_{i\in\set M}g^\pi\left({a_ix_i},{z_i}\right)\\
			+
			g^\pi\left( {\smashoperator{\sum_{i\in\set N_+}}a_ix_i\!-\!\smashoperator{\sum_{i\in \set L\cup\set M\cup\set N_-}}a_ix_i},{\sum_{i\in \set U}z_i} \right).
		\end{multlined}
	\end{align}
{If $\sum_{i\in \set N_+}a_ix_i>\sum_{i\in\set N_-}a_ix_i$ but such partition does not exist, then} $(t,x,z)\in\conv({\setRankOne})$ if and only if $t\ge f(x).$
	\item If  $\sum_{i\in \set N_+}a_ix_i\le\sum_{i\in\set N_-}a_ix_i$ and there exists an partition $\set L\cup\set M\cup\set U$ of $\set N_-\cap\support(x){\cap}\support(z)$ such that  
	\begin{align*}
		&1-\smashoperator{\sum_{i\in\set M\cup\set U}}z_i> 0,\quad\max_{i\in\set L}\frac{a_ix_i}{z_i}<\frac{\sum_{i\in \set L}a_ix_i}{1\!\!-\!\!\sum_{i\in\set M\cup\set U}z_i}\le\min_{i\in\set M}\frac{a_ix_i}{z_i}\\ 
		&\max_{i\in\set M}\frac{a_ix_i}{z_i}\le \frac{\sum_{i\in\set N_-}a_ix_i-\sum_{i\in \set L\cup\set M\cup\set N_+}a_ix_i}{\sum_{i\in\set U}z_i}<\min_{i\in\set U}\frac{a_ix_i}{z_i},
	\end{align*}
	then $(t,x,z)\in\conv({\setRankOne})$ if and only if 
	\begin{align*}
		\begin{multlined}[b][.8\textwidth]
			t\ge g^\pi\left( -{\sum_{i\in \set L}a_ix_i},{1\!-\!\smashoperator{\sum_{i\in\set M\cup\set U}}z_i} \right)+\sum_{i\in\set M}g^\pi\left(-{a_ix_i},{z_i}\right)\\+g^\pi\left( -{\smashoperator{\sum_{i\in\set N_-}}a_ix_i\!+\!\smashoperator{\sum_{i\in \set \set L\cup\set M\cup\set N_+}}a_ix_i},{\sum_{i\in \set U}z_i} \right).
		\end{multlined}
	\end{align*}
{If $\sum_{i\in \set N_+}a_ix_i\leq\sum_{i\in\set N_-}a_ix_i$ but such partition does not exist, then} $(t,x,z)\in\conv({\setRankOne})$ if and only if $t\ge f(x).$
	\end{itemize}
\end{theorem}
Note that $f(x)$ can be rewritten as $f(x)=\tilde g\left(\sum_{i\in \set N_-}a_ix_i-\sum_{i\in\set N_+}a_ix_i\right)$, where $\tilde g(s)=g(-s)$. Due to this symmetry, it suffices to prove the first assertion of Theorem~\ref{thm:explicit-pos}.
 Without loss of generality, we also assume $a_i=1\;\forall i\in [n]$; otherwise we can consider $\tilde x_i=a_ix_i$. 
 
Our proof proceeds as follows. First, we assume that the function $g$ satisfies additional regularity conditions, and that the point $(t,x,z){\in \R^{2n+1}}$ given in Theorem~\ref{thm:explicit-pos} is strictly positive. In other words, we provide a partial description of $\conv \set (\setRankOne)$ (incomplete since the boundary of the closure of the convex hull is not fully characterized). This first part of the proof is given in \S\ref{sec:partial}. Then in \S\ref{sec:fullDescription}, using approximations and limiting arguments, we provide a full description of $\conv \set (\setRankOne)$ for the general case. 
 
 \subsubsection{Partial characterization under regularity conditions}\label{sec:partial}
{We first formalize in Assumption~\ref{assum:r2} the regularity conditions we impose in this section.}
 \begin{assumption}\label{assum:r2}
 	In addition to Assumption~\ref{assum:r1}, function $g$ is a strongly convex and differentiable function with $g(0)=0=\min_{s\in\R}g(s)$. Moreover, $x>0$, $0<z\le \ones$ and $\sum_{i\in \set N_+}a_ix_i>\sum_{i\in\set N_-}a_ix_i$.
 \end{assumption}
{As mentioned above, the condition $\sum_{i\in \set N_+}a_ix_i>\sum_{i\in\set N_-}a_ix_i$ is without loss of generality due to symmetry. The other conditions in Assumption~\ref{assum:r2} are restrictive, but we show on \S\ref{sec:fullDescription} that they are not necessary. Finally, we point out that the point $(t,x,z)\in \R^{2n+1}$ given in Theorem~\ref{thm:explicit-pos} is not necessarily integral in $z$ --indeed, any points in the interior of $\conv ({\setRankOne})$ are not integral in $z$. In particular, the assumption $x_i>0$ does not imply that $z_i=1$.}

Observe that the positivity assumptions $x>0$ and $z>0$ imply that $\set N_+\cap\support(x)\cap\support(z)$ is simply $\set N_+$.
{Recall that }for any vector $y$ and index set {$\set J\subseteq [n]$, we denote $y(\set J)=\sum_{j\in\set J}y_j$ for convenience.} 

The workhorse of the derivation of Theorem~\ref{thm:explicit-pos} is the following minimization problem induced by the extended formulation \eqref{eq:rank-one} of $\conv({\setRankOne})$: for a given  $(x,z)$ such that $x> 0, 0<z\le \ones$ and  $\sum_{i\in \set N_+}x_i>\sum_{i\in\set N_-}x_i$,
\begin{equation}\label{eq:ext-pos}
	\begin{aligned}
		\min_{\lambda,\tau}\;&\sum_{i\in\set N_+}\lambda_ig\left( \frac{x_i-\tau_i}{\lambda_i} \right)+\sum_{i\in\set N_-}\lambda_ig\left(-\frac{x_i-\tau_i}{\lambda_i}\right)\\
		\text{s.t. }&\tau(\set N_+)-\tau(\set N_-)=0,\;
		0\le \tau\le x\\
		&\ones^\top \lambda\le1,\;
		0\le \lambda\le z.
	\end{aligned}
\end{equation}
{Indeed, a point $(t,x,z)$ satisfies the conditions in Proposition~\ref{prop:rank-one} if and only if $t$ is greater or equal than the optimal objective value of \eqref{eq:ext-pos}. Thus, projecting out variables $(\lambda,\tau)$ amounts to characterizing in closed form the optimal solutions of problem~\eqref{eq:ext-pos}.
	
	 Lemma~\ref{lem:pos-reduce-Nneg} below describes which bound constraints associated with variables $(\lambda,\tau)$ are tight (or not) in optimal solutions of \eqref{eq:ext-pos}.}
\begin{lemma}\label{lem:pos-reduce-Nneg}
	Under Assumption~\ref{assum:r2}, {there exists an optimal solution of $(\tau,\lambda)$ \eqref{eq:ext-pos} such that $\tau_i=x_i,\;\lambda_i=0\;\forall i\in\set N_-$ and $\lambda_i>0,\;\tau_i<x_i\;\forall i\in\set N_+$.}
\end{lemma}
\begin{proof}
	 For contradiction, we assume there exists some $i_-\in\set N_-$ such that $\tau_{i_-}<x_{i_-}$ in the optimal solution to \eqref{eq:ext-pos}. Since \[ \sum_{i\in\set N_+}\tau_i=\sum_{i\in\set N_-}\tau_i<\sum_{i\in\set N_-}x_i<\sum_{i\in\set N_+}x_i, \]
	there must be some $i_+\in\set N_+$ such that $\tau_{i_+}<x_{i_+}$. Since $g$ is convex and attains its minimum at $0$, $g$ is increasing over $[0,+\infty)$ and decreasing over $(-\infty, 0]$. It follows that increasing $\tau_{i_-}$ and $\tau_{i_+}$ by the same sufficiently small amount of $\epsilon>0$ would improve the objective function, which contradicts with the optimality. Hence, $\tau_i=x_i\;\forall i\in \set N_-$. It follows that one can safely take $\lambda_i=0\;\forall i\in\set N_-$.
	
	We now prove the second part of the statement, namely, that there exists an optimal solution with $\lambda>0$. If $\lambda_i=0,i\in\set N_+$, since $g$ is strongly convex, one must have $\tau_i=x_i$, otherwise $g^\pi(x_i-\tau_i,0)=+\infty$. Moreover, if $\tau_i=x_i$, one can safely take $\lambda_i=0$. Finally, assume $\lambda_i=0,\tau_i=x_i$ for some index $i$. In this case, constraints $\tau(\set N_+)-\tau(\set N_-)=0$, the previously proven property that $\tau(\set N_-)=x(\set N_-)$ and the assumption $x(\set N_+)>x(\set N_-)$ imply that there exists an index $j$ where $\tau_j<x_j$ and $\lambda_j>0$. Setting $(\tilde \lambda_i,x_i-\tilde \tau_i)=\epsilon(\lambda_j,x_j-\tau_j)$ and $(\tilde \lambda_j,x_j-\tilde \tau_j)=(1-\epsilon)(\lambda_j,x_j-\tau_j)$ for some small enough $\epsilon>0$, we find that the new feasible solution is still optimal since $g^\pi(\cdot,\cdot)$ is a homogeneous function. The conclusion follows.
\end{proof}
By changing variable $\tau\leftarrow x-\tau$, it follows from Lemma~\ref{lem:pos-reduce-Nneg} that \eqref{eq:ext-pos} can be simplified to 
\begin{subequations}\label{eq:ext-pos-reduced}
	\begin{align}
		\min_{\lambda,\tau}\;&\sum_{i\in\set N_+}\lambda_ig\left( \frac{\tau_i}{\lambda_i} \right)\\
		\text{s.t. }&\tau(\set N_+)=C\\
		&0< \tau\le x\\
		&\lambda(\set N_+)\le1\\
		&0< \lambda\le z,
	\end{align}
\end{subequations}
where $C\defeq x(\set N_+)-x(\set N_-)$. Denote the derivative of $g(t)$ by $g'(t)$. Define $G(t)\defeq tg'(t)-g(t)>0$ and $\hat g(t)\defeq tg(1/t)$. {To characterize the optimal solutions of \eqref{eq:ext-pos-reduced}, we will explicitly solve the KKT conditions associated with the problem. Such conditions include equations involving $g'$ and $G$.}
The next lemma shows that these functions are monotone, which is essential to solve the KKT system in closed form.
\begin{lemma}\label{lem:explicit-pos-monotonicity}
	Under Assumption~\ref{assum:r2}, $g'(t)$ and $G(t)$ are increasing over $(0,+\infty)$ and thus invertible. Moreover, $\hat g(t)$ is strictly convex over $(0,+\infty)$.
\end{lemma}
\begin{proof}
	Strict monotonicity of $g'(t)$ over $[0,+\infty)$ follows directly from the strict convexity of $g$. Since $g$ is strongly convex, it can be written as $g(t)=g_1(t)+g_2(t)$, where $g_1(t)=\alpha t^2$ is quadratic with a certain $\alpha>0$ and $g_2(t)$ is convex.
	For any $t>0$,
	\begin{align*}
		&G(t+\epsilon)-G(t)=\epsilon g'(t+\epsilon)-(g(t+\epsilon)-g(t))+t\left( g'(t+\epsilon)-g'(t) \right)\\
		=&\epsilon\left( g'(t+\epsilon)-g'(t) \right)+t\left( g'(t+\epsilon)-g'(t) \right) +\littleo(\epsilon)\tag{Taylor expansion}\\
		=&\epsilon\left(g'(t+\epsilon)-g'(t)\right)+t(g_2'(t+\epsilon)-g_2'(t))+2\alpha t\epsilon+\littleo(\epsilon)\\
		\ge&2\alpha t\epsilon+\littleo(\epsilon)>0\tag{Monotonicity of the derivative of a convex function},
	\end{align*}
as $\epsilon$ is small enough. Namely, $G(t)$ is  an increasing function over $[0,+\infty)$. Moreover, $G(0)=0$ implies that $G(t)>0\;\forall t\in(0,+\infty)$. To prove the last conclusion, $\hat g'(t)=g(1/t)-g'(1/t)/t=-G(1/t)$ is increasing over $(0,+\infty)$, which implies the strict convexity of $\hat g$.
\end{proof}
Assume $(\lambda,\tau)$ is an optimal solution to \eqref{eq:ext-pos-reduced}. Let $r_i\defeq\frac{\tau_i}{\lambda_i}\;\forall i\in\set N_+$. The next lemma reveals the structure of the optimal solution to \eqref{eq:ext-pos-reduced}. Namely, unless {the} $r_i$'s are identical, either $\lambda_i$ or $\tau_i$ attains the upper bound. 
\begin{lemma}\label{lem:explicit-pos-alternative}
	Under Assumption~\ref{assum:r2}, if $r_i>r_j$, then $\lambda_i=z_i$ and  $\tau_j=x_j$ in the optimal solution to \eqref{eq:ext-pos-reduced}.
\end{lemma}
\begin{proof}
	For contradiction we assume $\lambda_i<z_i$. Take $\epsilon>0$ small enough and let $\epsilon_i=\epsilon/\tau_i$ and $\epsilon_j=\epsilon/\tau_j$. Since $\hat g$ is strictly convex and $1/r_i<1/r_j$, one can deduce that
	\begin{align*}
		&\frac{\hat g(1/r_i+\epsilon_i)-\hat g(1/r_i)}{\epsilon_i}<\frac{\hat g(1/r_j)-\hat g(1/r_j-\epsilon_j)}{\epsilon_j}\\
		\Leftrightarrow\quad&\begin{aligned}
			\begin{multlined}[b][.87\textwidth]
				\frac{(\lambda_i+\epsilon)/\tau_ig(\tau_i/(\lambda_i+\epsilon))-\lambda_i/\tau_ig(\tau_i/\lambda_i)}{\epsilon/\tau_i}<\\\frac{\lambda_j/\tau_jg(\tau_j/\lambda_j)-(\lambda_j-\epsilon)/\tau_jg(\tau_j/(\lambda_j-\epsilon))}{\epsilon/\tau_j}
			\end{multlined}
		\end{aligned}\\
		\Leftrightarrow\quad&(\lambda_i+\epsilon)g\left( \frac{\tau_i}{\lambda_i+\epsilon}\right)+(\lambda_j-\epsilon)g\left(\frac{\tau_j}{\lambda_j-\epsilon}\right)-\lambda_ig\left(\frac{\tau_i}{\lambda_i}\right)-\lambda_jg\left(\frac{\tau_j}{\lambda_j}\right)<0.
	\end{align*}
It implies that we can improve the objective value of \eqref{eq:ext-pos-reduced} by increasing $\lambda_i$ and decreasing $\lambda_j$ by $\epsilon$, which contradicts with the optimality. Hence, one can deduce $\lambda_i=z_i$.

Similarly, the second conclusion follows by 
\begin{align*}
	&\lambda_ig\left(\frac{\tau_i-\epsilon}{\lambda_i}\right)+\lambda_jg\left(\frac{\tau_j+\epsilon}{\lambda_j}\right)-\lambda_ig\left( \frac{\tau_i}{\lambda_i} \right)-\lambda_jg\left(\frac{\tau_j}{\lambda_j}\right)\\
	=&(g'(r_j)-g'(r_i))\epsilon+\littleo(\epsilon)<0,
\end{align*}
since $g'(\cdot)$ is strictly increasing.
\end{proof}

We are now ready to prove Theorem~\ref{thm:explicit-pos} {under the additional regularity conditions.}
\begin{proposition}\label{prop:explicit-pos-merit}
	Under Assumption~\ref{assum:r2}, the conclusion in Theorem~\ref{thm:explicit-pos} holds true.
\end{proposition}
\begin{proof}
	We discuss two cases defined by $z(\set N_+)$ separately.
	
\noindent\emph{Case 1: $z(\set N_+)> 1$.} Since the objective function of \eqref{eq:ext-pos-reduced} is decreasing with respect to $\lambda$ by Lemma~\ref{lem:monotonicity-persp}, one must have $\lambda(\set N_+)=1$. Thus, program \eqref{eq:ext-pos-reduced} can be reduced to
\begin{subequations}\label{eq:ext-pos-case1}
	\begin{align}
		\min_{\lambda,\tau}\;&\sum_{i\in\set N_+}\lambda_ig\left( \frac{\tau_i}{\lambda_i} \right)\label{eq:ext-pos-case1-obj}\\
		\text{s.t. }&\tau(\set N_+)=C\tag{$\alpha$}\\
		&0< \tau\le x\tag{$\beta$}\\
		&\lambda(\set N_+)=1\tag{$\delta$}\\
		&0< \lambda\le z\tag{$\gamma$}.
	\end{align}
\end{subequations}
Assume $(\tau,\lambda)$ is the optimal solution to \eqref{eq:ext-pos-case1} and define $r_i=\tau_i/\lambda_i$. There are two possibilities -- either all $r_i$'s are identical or there are at least two distinct values of $r_i$'s. In the former case, denote $r=r_i\;\forall i\in\set N_+$, that is, $\tau_i=r\lambda_i\;\forall i\in\set N_+$. Then $\tau(\set N_+)/\lambda(\set N_+)=r\lambda(\set N_+)/\lambda(\set N_+)=r$, and in particular $r=C$. In this case, \eqref{eq:ext-pos-case1-obj} reduces to $\lambda(\set N_+)g(r)=g(C)$.

Now we assume there are at least two distinct values of $r_i$'s. By Lemma~\ref{lem:explicit-pos-alternative}, for all $i\in\set N_+$, either $x_i=\tau_i$ or $\lambda_i=z_i$. It follows that $\set N_+=\set L\cup\set M\cup\set U$, where 
\begin{equation}\label{eq:explicit-pos-partition}
	\set L=\{ i:\tau_i=x_i,\lambda_i<z_i \},\;\set M=\{i:\tau_i=x_i,\lambda_i=z_i  \},\;\set U=\{ \tau_i<x_i,\lambda_i=z_i \}.
\end{equation}
Since all constraints of \eqref{eq:ext-pos-case1} are linear, KKT conditions are necessary and sufficient for optimality of $(\lambda,\tau)$. Let $(\alpha,\beta,\gamma,\delta)$ be the dual variables associated with each constraint of \eqref{eq:ext-pos-case1}.  It follows that $\gamma_i=0\;\forall i\in\set L$ and $\beta_i=0\;\forall i\in\set U$. The KKT conditions of \eqref{eq:ext-pos-case1} can be stated as follows (left column: statement of the KKT condition; right column: equivalent simplification)
\begin{equation*}
	\begin{array}{l|lr}
		g'\left( \frac{x_i}{\lambda_i} \right)-\alpha+\beta_i=0&\beta_i=\alpha-g'\left( \frac{x_i}{\lambda_i} \right) &\forall i\in\set L\\
		g'\left(\frac{x_i}{z_i}\right)-\alpha+\beta_i=0&\beta_i=\alpha-g'\left( \frac{x_i}{z_i} \right)&\forall i\in\set M\\
		g'\left( \frac{\tau_i}{z_i} \right)-\alpha=0&\tau_i=z_i(g')^{-1}(\alpha)&\forall i\in\set U\\
		g\left( \frac{x_i}{\lambda_i}\right)-\frac{x_i}{\lambda_i}g'\left(\frac{x_i}{\lambda_i} \right)+\delta=0&\lambda_i=\frac{x_i}{G^{-1}(\delta)}&\forall i\in\set L\\
		g\left( \frac{x_i}{z_i}\right)-\frac{x_i}{z_i}g'\left(\frac{x_i}{z_i} \right)+\delta+\gamma_i=0&\gamma_i=G\left( \frac{x_i}{z_i} \right)-\delta&\forall i\in\set M\\
		g'\left(\frac{\tau_i}{z_i} \right)-\frac{\tau_i}{z_i}g'\left(\frac{\tau_i}{z_i} \right)+\delta+\gamma_i=0&\gamma_i=G\left( \frac{\tau_i}{z_i} \right)-\delta&\forall i\in\set U\\
		\tau(\set U)+x(\set M)+x(\set L)=C&\tau(\set U)=C-x(\set M)-x(\set L)\\
		\lambda(\set L)+z(\set M)+z(\set U)=1&\lambda(\set L)=1-z(\set M)-z(\set U)\\
		\multicolumn{3}{l}{0<\tau_i<x_i\;\forall i\in\set U,\;0<\lambda_i<z_i\;\forall i\in\set L,\;\beta_i\ge0\;\forall i\in\set L\cup\set M,\;\gamma_i\ge0\;\forall i\in\set M\cup\set U}.
	\end{array}
\end{equation*}
Denote by $\bar C=C-x(\set M)-x(\set L)$ and $\bar z=1-z(\set M)-z(\set U)$. Then
\begin{equation*}
	\begin{array}{lcl}
		\tau(\set U)=z(\set U)(g')^{-1}(\alpha)=\bar C&\Rightarrow& \alpha=g'\left( {\bar C}/{z(\set U)} \right)\\
		\lambda(\set L)=x(\set L)/G^{-1}(\delta)=\bar z&\Rightarrow& \delta=G\left( {x(\set L)}/{\bar z} \right).
	\end{array}
\end{equation*}
Thus, one can first substitute out $\alpha$ and $\delta$ to get $\beta_i\;\forall i\in\set M,\tau_i\;\forall i\in\set U,\lambda_i\;\forall i\in\set L,\gamma_i\;\forall i\in\set M$. Then one can plug in $\lambda_i\;\forall i\in\set L$ and $\tau_i\;\forall i\in\set U$ to work out $\beta_i\;\forall i\in\set L$ and $\gamma_i\;\forall i\in\set U$. Hence, we deduce that the KKT system is equivalent to 
\begin{subequations}
	\begin{align}
		&\beta_i=g'(\bar C/z(\set U))-g'(x(\set L)/\bar z)\ge0&\forall i\in\set L\label{eq:ineq-ul1}\\
		&\beta_i=g'(\bar C/z(\set U))-g'(x_i/z_i)\ge0&\forall i\in\set M\label{eq:ineq-um}\\
		&\tau_i=\bar Cz_i/z(\set U)\in(0,x_i)&\forall i\in\set U\label{eq:ineq-uu}\\
		&\lambda_i=\bar zx_i/x(\set L)\in(0,z_i)&\forall i\in\set L\label{eq:ineq-ll}\\
		&\gamma_i=G(x_i/z_i)-G(x(\set L)/\bar z)\ge 0&\forall i\in\set M\label{eq:ineq-ml}\\
		&\gamma_i=G( \bar C/z(\set U))-G(x(\set L)/\bar z)\ge 0&\forall i\in\set U\label{eq:ineq-ul2}.
	\end{align}
\end{subequations}

Because $g'$ and $G$ are increasing from Lemma~\ref{lem:explicit-pos-monotonicity}, the KKT system has a solution if and only if 
\begin{align*}
	\min_{i\in\set U}\frac{x_i}{z_i}\stackrel{\eqref{eq:ineq-uu}}{>}\frac{\bar C}{z(\set U)}\stackrel{\eqref{eq:ineq-um}}{\ge} \max_{i\in\set M}\frac{x_i}{z_i}\ge \min_{i\in\set M}\frac{x_i}{z_i}\stackrel{\eqref{eq:ineq-ml}}{\ge}\frac{x(\set L)}{\bar z}\stackrel{\eqref{eq:ineq-ll}}{>}\max_{i\in \set L}\frac{x_i}{z_i},
\end{align*}
which implies $\frac{\bar C}{z(\set U)}\ge\frac{x(\set L)}{\bar z}\Leftrightarrow$\eqref{eq:ineq-ul1} and \eqref{eq:ineq-ul2}.
Moreover, by using the solution to the KKT system, the optimal value of \eqref{eq:ext-pos-case1} is
\[ g\left( \frac{x(\set L)}{\bar z} \right)+\sum_{i\in\set M}z_ig\left( \frac{x_i}{z_i} \right)+z(\set U)g\left(\frac{\bar C}{z(\set U)}\right).\]
\noindent\emph{Case 2: $z(\set N_+)\le 1$.} Since the objective function of \eqref{eq:ext-pos-reduced} is decreasing with respect to $\lambda$, one must have $\lambda_i=z_i\;\forall i\in\set N_+$. Thus, problem~\eqref{eq:ext-pos-reduced} reduces to
\begin{subequations}\label{eq:ext-pos-case2}
	\begin{align}
		\min_{\tau}\;&\sum_{i\in\set N_+}z_ig\left( \frac{\tau_i}{z_i} \right)\label{eq:ext-pos-case2-obj}\\
		\text{s.t. }&\tau(\set N_+)=C\tag{$\alpha$}\\
		&0< \tau\le x\tag{$\beta$}.
	\end{align}
\end{subequations}
Assume $\tau$ is the optimal solution to \eqref{eq:ext-pos-case2}. Similarly, define 
\[ \set L=\emptyset,\quad\set M=\{i:\tau_i=x_i \},\quad\set U=\{ i:\tau_i<x_i \}. \]
Then the KKT conditions can be written as
\begin{equation*}
	\begin{array}{l|lr}
		g'\left( \frac{x_i}{z_i} \right)-\alpha+\beta_i=0&\beta_i=-g'\left( \frac{x_i}{z_i} \right)+\alpha&\forall i\in \set M\\
		g'\left(\frac{\tau_i}{z_i}\right)-\alpha&\tau_i=(g')^{-1}(\alpha)z_i&\forall i\in \set U\\
		x(\set M)+\tau(\set U)=C&\tau(\set U)=C-x(\set M)=\bar C=(g')^{-1}(\alpha)z(\set U)\\
		\multicolumn{3}{l}{0<\tau_i<x_i\;\forall i\in \set U,\;\beta_i\ge0\;\forall i\in\set M.}
	\end{array}
\end{equation*}
It follows that $\alpha=g'(\bar C/z(\set U))$. Plugging $\alpha$ in, one arrives at
\begin{align*}
	&\beta_i=g'(\bar C/z(\set U))-g'(x_i/z_i)\ge0&\forall i\in\set M\\
	&\tau_i =z_i\bar C/z(\set U)\in(0,x_i)&\forall i\in \set U.
\end{align*}
Therefore, the KKT system has a solution if and only if 
\[ \min_{i\in\set U}\frac{x_i}{z_i}>\frac{\bar C}{z(\set U)}\ge\max_{i\in\set M}\frac{x_i}{z_i}. \]
The proof is finished.
\end{proof}

\subsubsection{Full description of the closure of the convex hull}\label{sec:fullDescription}
From Proposition~\ref{prop:explicit-pos-merit}, we see that Theorem~\ref{thm:explicit-pos} holds if additional {regularity conditions} are imposed -- namely $g$ is finite, strongly convex and differentiable, $0=g(0)=\min_{t\in\R}g(t)$ and $x>0,z>0$. To complete the proof of the theorem, we now show how to remove the assumptions, one by one. {In particular, we show how to approximate an arbitrary function $g$ with a function that satisfies the regularity conditions without changing the optimal solutions $(\lambda^*,\tau^*)$ of the optimization problems \eqref{eq:ext-pos}. Similarly, if the assumption  $z>0$ is not satisfied, then we prove the result by constructing a sequence of positive solutions converging to $(t,x,z)$. } 

\begin{proof}[Proof of Theorem~\ref{thm:explicit-pos}]
	 Due to the symmetry mentioned above, we only prove the first conclusion in the theorem under Assumption~\ref{assum:r1} and $a=\ones$.

	 A key observation is that the optimal primal solution to \eqref{eq:ext-pos-reduced}, given by \eqref{eq:ineq-uu} and \eqref{eq:ineq-ll}, does not involve function $g$ and only relies on $x$ and $z$ (while the values of the dual variables does depend on $g$). Denote this optimal solution by $(\tau^*,\lambda^*)$ and the objective function of \eqref{eq:ext-pos-reduced} by $h(\tau,\lambda;g)$. Then $h(\tau^*,\lambda^*;g)\le h(\tau,\lambda;g)$ for all feasible solutions $(\tau,\lambda)$ to \eqref{eq:ext-pos-reduced} and all functions $g$ satisfying Assumption~\ref{assum:r2}.
	
	{First suppose that $g$} is not a strongly convex function. {We show that the solution $(\tau^*,\lambda^*)$ above is still optimal for \eqref{eq:ext-pos-reduced}.}
 Consider $ g_\epsilon(s)=g(s)+\epsilon s^2$, where $\epsilon>0$. Since $ g_\epsilon$ is strongly convex, { $(\tau^*,\lambda^*)$ is optimal if $g$ is replaced with $g_\epsilon$. Optimality of $(\tau^*,\lambda^*)$ is equivalent (by definition)} to 
	$ h(\tau^*,\lambda^*;g_\epsilon)\le h(\tau,\lambda;g_\epsilon)$ for any feasible solution $(\tau,\lambda)$ of \eqref{eq:ext-pos-reduced}. {Letting $\epsilon\to 0$, we find that}
 $h(\tau^*,\lambda^*;g)\le h(\tau,\lambda;g)${, and in particular $(\tau^*,\lambda^*)$ is optimal under function $g$}. {Thus, Proposition~\ref{prop:explicit-pos-merit} holds even} if $g$ is a {(not necessarily strongly convex)} differentiable function with $0=g(0)=\min_{t\in\R}g(t)$.
	
	{Now suppose that $g$} is not a differentiable function. Consider its \emph{Moreau-Yosida regularization} $e_\epsilon g(s)\defeq \min_w\left\{g(w)+\frac{1}{2\epsilon}(s-w)^2 \right\}\le g(s),$ where $\epsilon>0${, which} is a differentiable convex function; see Corollary~4.5.5,~ \cite{hiriart2004fundamentals}. Moreover, $e_\epsilon g(s)\to g(s)$ as $\epsilon\to 0$; see Theorem~1.25, \cite{rockafellar2009variational}. {Similarly to the previous argument, }since  $ h(\tau^*,\lambda^*;e_\epsilon g)\le h(\tau,\lambda;e_\epsilon g)$ {for any feasible solution $(\tau,\lambda)$ of \eqref{eq:ext-pos-reduced}}, letting $\epsilon\to 0${, one can deduce that $h(\tau^*,\lambda^*)$ is optimal under objective $g$.} {Thus, Proposition~\ref{prop:explicit-pos-merit} holds even} if $g$ is a  non-differentiable finite convex function with $0=g(0)=\min_{t\in\R}g(t)$.
	
	For a general finite convex function $g$ with $g(0)=0$, $\tilde g(s)=g(s)-cs$ is a convex function with $0=\tilde g(0)=\min_{s\in\R}\tilde g(s)$, where $c\in\partial g(0)$. The conclusion follows by applying the theorem to $\tilde{g}$.

	{To remove the positivity assumptions, first note that} if there exists $i\in\set N_+$ such that $x_i=0$, then $\tau_i=0$ in \eqref{eq:ext-pos} which implies that one can safely set $\lambda_i=0$ in \eqref{eq:ext-pos}. Hence, we can exclude the variables associated with index $i$ from consideration and reduce the problem to a lower-dimensional case. For this reason, without loss of generality, we assume $x_i>0\;\forall i\in[n]$. 
	
	{Finally,} define $\set N_0\defeq\left\{i\in\set N_+:z_i=0\right\}$, $\tilde f(x,z)$ as the RHS of \eqref{eq:convex-envolope}, and $f^\dagger(x,z)$ as the optimal value of \eqref{eq:ext-pos}. Let $r>0$ be a sufficiently large number and consider $z^r$ defined as $z^r_i=x_i/r$ if $i\in\set N_0$ and $z^r_i=z_i$ otherwise. Note that $\lim\limits_{r\to+\infty}z^r=z$. If there exists a partition $\set L\cup\set M\cup\set U$ of $\set N_+\cap\support(z)$ stated in the theorem, then $\set L\cup\set M\cup\tilde{\set U}$ is the partition of $\set N_+$ associated with $(x,z^r)$ where $\tilde{\set U}=\set U\cup\set N_0$. Since $z^r>0$, the conclusion holds for $(x,z^r)$, i.e.  $f^\dagger\left(x,z^r\right)=\tilde f\left(x,z^r\right)$. Because $f^\dagger$ and $\tilde f$ are closed convex functions, $f^\dagger(x,z)=\lim\limits_{r\to+\infty} f^\dagger\left(x,z^r\right)=\lim\limits_{r\to+\infty} \tilde f\left(x,z^r\right)=\tilde f\left(x,z\right).$ On the other hand, if for $(x,z^r)$ such a partition $\set L\cup\set M\cup\set U$ of $\set N_+$ does not exist, then neither does for $(x,z)$ because otherwise, $\left(\set L\backslash \set N_0,\set M\backslash \set N_0,\set U\backslash \set N_0\right)$ would be a proper partition of $\set N_+\cap\support(z)$. In this case, the conclusion follows from the closedness of $f$ and $f^\dagger$. This completes the proof.
	
\end{proof}
\begin{remark}
Sets $\set L$, $\set M$ and $\set U$ in Theorem~\ref{thm:explicit-pos} can be found in $\bigO{n^2}$ time. Indeed, without loss of generality, we assume that  $\sum_{i\in \set N_+}a_ix_i>\sum_{i\in\set N_-}a_ix_i$. First, sort and index $x_i/z_i$ in a nondecreasing order. It follows from the conditions in Theorem~\ref{thm:explicit-pos} that if such $\set L$, $\set M$ and $\set U$ exist, then there must be some $k_1$ and $k_2$ such that $\set L=\{ i\in\set N_+:i<k_1 \}$, $\set M=\{ i\in\set N_+:k_1\le i\le k_2 \}$ and $\set U=\{ i\in\set N_+:i>k_2 \}$. Consequently, one can verify the conditions in Theorem~\ref{thm:explicit-pos} by enumerating all possible combinations $\{k_1,k_2\}\subseteq \set N_+$.
\end{remark}
Now we turn to the special case where $\set N_-=\emptyset$, that is, every entry of $a$ is positive.
\begin{corollary}\label{cor:explicit-pp}
	Under Assumption~\ref{assum:r1} and $a>0$, point  $(t,x,z)\in\conv({\setRankOne})$ if and only if $(x,z)\in[0,1]^n\times\R_+^n$ and there exists a partition $\set L\cup\set M=\support(x)\cap\support(z)$ such that
	\begin{equation}\label{eq:explicit-pp-condition}
		1-\smashoperator{\sum_{i\in\set M}}z_i\ge0,\quad\max_{i\in\set L}\frac{a_ix_i}{z_i}<\frac{\sum_{i\in\set L}a_ix_i}{1-\sum_{i\in\set M}z_i}\le \min_{i\in\set M}\frac{a_ix_i}{z_i},
	\end{equation}
	and the following inequality holds
	\begin{equation}\label{eq:explcit-pp-epi}
		t\ge g^\pi\left(\sum_{i\in\set L}a_ix_i,1-\sum_{i\in\set M}z_i \right)+\sum_{i\in \set M}g^\pi(a_ix_i,z_i)+g^\pi\left(\smashoperator{\sum_{i\in\set N_+}}a_ix_i-\smashoperator{\sum_{i\in\set L\cup\set M}}a_ix_i,0\right).
	\end{equation}
\end{corollary}
\begin{proof}
	In this particular setting, it is easy to see that $\set N_-=\emptyset$ and $\tau_i=0\;\forall i\in[n]$ in \eqref{eq:ext-pos}. Thus, the partition defined in \eqref{eq:explicit-pos-partition} always exists with $\set U=\emptyset$. The conclusion follows from Theorem~\ref{thm:explicit-pos}. 
\end{proof} 

Finally, we close this section by generalizing the main result of \cite{atamturk2018strong} to non-quadratic functions. Specifically, in \cite{atamturk2018strong}, the authors studied the set 
\[ {\setRankOne^2}\defeq\left\{ (t,x,z)\in\R^{3}\times{\{0,1\}^n}:\quad\begin{aligned}
&t\ge g(a_1x_1-a_2x_2), x_i\ge 0,\;i=1,2,\\
&x_i(1-z_i)=0,\;i=1,2
\end{aligned}\right\}, \]
where $g$ is quadratic, and provided the description of $\conv({\setRankOne^2})$ in the original space of variable. A similar result holds for general convex functions.
\begin{corollary}
	Given a convex function $g(\cdot)$ with $\dom(g)=\R$ and $f(x_1,x_2)=g\left(a_1x_1-a_2x_2\right)$, where $a>0$, point  $(t,x,z)\in\conv({\setRankOne^2})$ if and only if $(x,z)\in[0,1]^2\times\R_+^2$ and
	\begin{align*}
		t\ge\begin{cases}
			g^\pi(a_1x_1-a_2x_2,z_1)&\text{if }a_1x_1\ge a_2x_2\\
			g^\pi(a_2x_2-a_1x_1,z_2)&\text{if }a_2x_2\ge a_1x_1.
		\end{cases}
	\end{align*}
\end{corollary}
\begin{proof}
	In this case, $\set N_+$ is a singleton in \eqref{eq:ext-pos-reduced}.
\end{proof}

\subsection{Implementation}\label{sec:implementation}

In this section, we discuss the implementation of the results given in Theorem~\ref{thm:explicit-pos} (for the quadratic case) with conic quadratic solvers. 

 A key difficulty towards using the convexification is that inequalities \eqref{eq:convex-envolope} are not valid: while they describe $\conv(Q)$ in their corresponding region, determined by partition $\set L\cup \set M\cup \set U$, they may cut off points of $\conv(Q)$ elsewhere. To circumvent this issue, \citet{atamturk2020supermodularity} propose valid inequalities, each requiring $\mathcal{O}(n)$ additional variables and corresponding exactly with \eqref{eq:convex-envolope} in the corresponding region, and valid elsewhere. The inequalities are then implemented as cutting surfaces, added on the fly as needed. It is worth nothing that since the optimization problems considered are nonlinear, and convex relaxations are solved via interior point solvers, adding a cut requires resolving again the convex relaxation (without the warm-starting capabilities of the simplex method for linear optimization). 

In contrast, we can use Proposition~\ref{prop:rank-one} directly to implement the inequalities. When specialized to quadratic functions $g$, and with the introduction of auxiliary variables $u$ to model conic quadratic cones, we can restate Proposition~\ref{prop:rank-one} as: $(x,y,t)\in \conv (\set Q)$ if and only if there exists $(\lambda,\tau,u)\in \R^{3n}$ such that 
	\begin{subequations}\label{eq:conicQuadraticRepresentability}
	\begin{align}
	&t\ge \sum_{i=1}^na_i^2u_i,\\
	&\lambda_i u_i\geq (x_i-\tau_i)^2,\; u_i\geq 0,\;\forall i\in [n],\label{eq:conicQuadraticRepresentability_rotated}\\ 
	&a^\top \tau=0,\;0\le \tau_i\le x_i\;\forall i\in\set I_+,\\
	&\lambda_i\le z_i\le 1,\;\forall i\in[n],\\
	&\lambda\ge 0,\;\sum_{i=1}^n \lambda_i\le1.
	\end{align}
	\end{subequations}
	is satisfied. Inequalities \eqref{eq:conicQuadraticRepresentability_rotated} are (convex) rotated cone constraints, which can be handled by most off-the-shelf conic quadratic solvers, and every other constraint is linear. Note that using \eqref{eq:conicQuadraticRepresentability} requires adding $\mathcal{O}(n)$ variables \emph{once} --instead of adding a similar number of variables \emph{per inequality added}, with exponentially many inequalities required to describe $\conv(\set Q)$--, and thus is a substantially more compact formulation than the one presented in \cite{atamturk2020supermodularity}. 
	
\begin{figure}
		\includegraphics[trim={10cm 6cm 10cm 6cm}, clip,width=\textwidth]{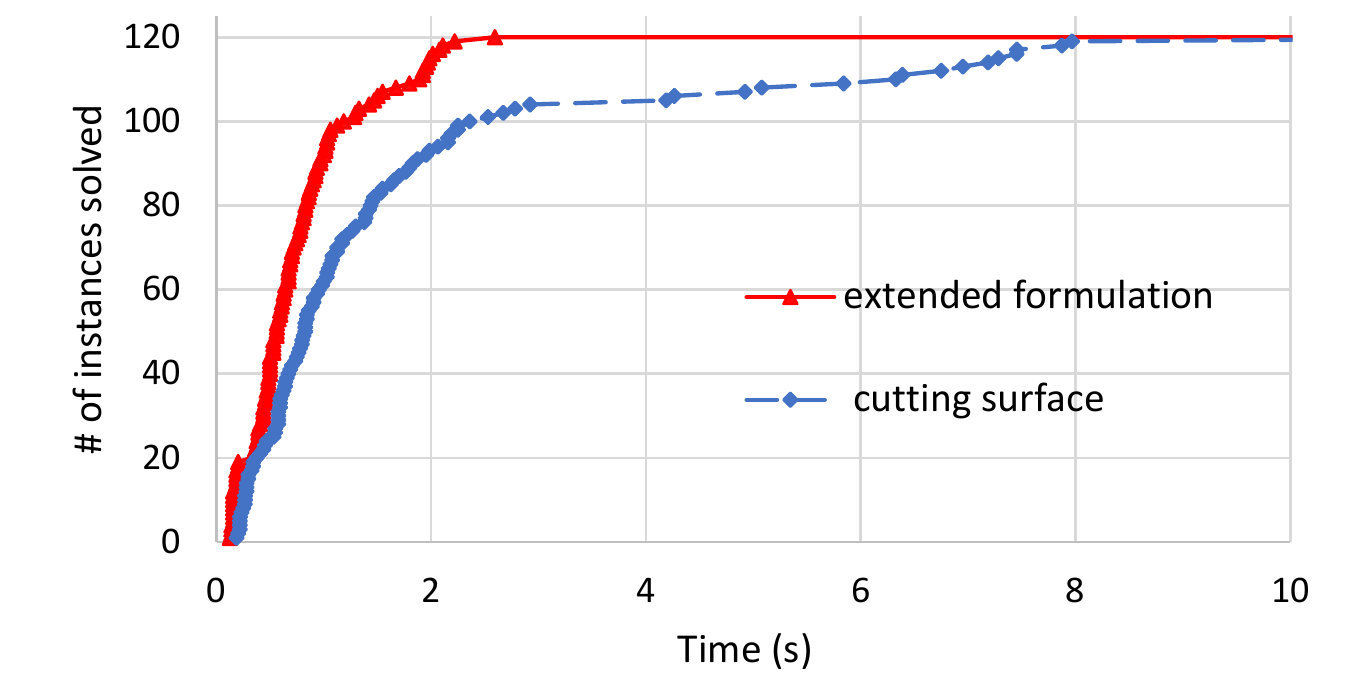}
		\caption{Number of instances solved as a function of time. The cutting surface requires on average 1.79 seconds to solve an instance, and can solve all 120 instances in 13.3 seconds or less. In contrast, the extended formulation \eqref{eq:conicQuadraticRepresentability} requires on average 0.78 seconds, and solves all instances in less than 2.6 seconds.}
		\label{fig:performance profile}
	\end{figure}

	To illustrate the benefits resulting from a more compact formulation, we compare the two formulations in instances used by \cite{atamturk2020supermodularity}, available online at \url{https://sites.google.com/usc.edu/gomez/data}. The instances correspond to portfolio optimization problems of the form 
	\begin{subequations}\label{eq:portfolio}
	\begin{align}
	\min\;&\sum_{k=1}^K t_k+\sum_{i=1}^n (d_ix_i)^2 \label{eq:portfolio_obj}\\
	\text{s.t.}\;&t_k\geq (a_k^\top x)^2&{\forall k\in [K]}\label{eq:portfolio_rankone}\\
	&\ones^\top x=1\\
	&c^\top x- h^\top z\geq b\\
	&0\leq x\leq z\\
	&x\in \R^n,\; z\in \{0,1\}^n,
	\end{align} 
	\end{subequations}
	where $a_k\in \R^n$ for all $k\in [K]$, $c,d,h\in \R_+^n$ and $b\in \R_+$. Strong relaxations can be obtained by relaxing the integrality constraints  $z\in \{0,1\}^n\to z\in [0,1]^n$, using the perspective reformulation {by replacing $(d_ix_i)^2$ with $ (d_ix_i)^2/z_i$} in the objective {\eqref{eq:portfolio_obj}}, and adding inequalities \eqref{eq:conicQuadraticRepresentability} (or using cutting surfaces) corresponding to each rank-one constraint \eqref{eq:portfolio_rankone}. Figure~\ref{fig:performance profile} summarizes the computational times required to solve the convex relaxations across 120 instances with $n\in \{200,500\}$ and $K\in \{5,10,20\}$, using the CPLEX solver in a laptop with Intel Core i7-8550U CPU and 16GB memory. In short, the extended formulation \eqref{eq:conicQuadraticRepresentability} is on average twice as fast as the cutting surface method proposed in \cite{atamturk2020supermodularity}, and up to five times faster in the more difficult instances (in addition to arguably being easier to implement). 
	
We also tested the effect of the extended formulation using CPLEX branch-and-bound solver. While we did not encounter numerical issues resulting in incorrect behavior by the solver (the cutting surface method does result in numerical issues, see \cite{atamturk2020supermodularity}), the performance of the branch-and-bound method is substantially impaired when using the extended formulation. We discuss in more detail the issues of using the extended formulation with CPLEX branch-and-bound method in the appendix.

\section{\np-hardness with bound constraints on continuous variables}\label{sec:complexity}
The set $\set Q$ {and its rank-one specialization $\set Q_1$} studied so far assume that the continuous variables are either unbounded, or non-negative/non-positive. Either way, set $\set Q$ {and $\set Q_1$} admit a similar disjunctive form given in Theorem~\ref{thm:main} and Proposition~\ref{prop:rank-one}{, respectively}. A natural question is whether the addition of bounds on the continuous variables results in similar convexifications, or if the resulting sets are structurally different. In this section, we show that {it would be unlikely} to {find a tractable description of the} the convex hulls with bounded variables, unless $\set P=\np$.

Consider the set
\[ \set Q_B\defeq\left\{ (t,x,z)\in\R^{n+1}\times\{0,1\}^n: t\ge f(x),\;0\le x\le z \right\}. \]
We show that {optimizing over set} $\conv(\set Q_B)$ is \np-hard even when $\rank(f)=1$. Two examples are given to illustrate this point -- \emph{single node flow sets} and {rank-one quadratic forms}.
\subsection*{Single-node fixed-charge flow set} The single-node fixed-charge flow set is the mixed integer linear set defined as 
\[ \set T\defeq\left\{ (x,z)\in\R^n\times\{0,1\}^n:\sum_{i=1}^{n}a_ix_i\le b,\; 0\le x\le z \right\}, \]
where $0<a_i\le b\;\forall i\in[n]$. Note that one face of the single-node flow set $\text{conv}(\set T\cap\{(x,z):x_i=z_i\;\forall i\in[n]\})$ is isomorphic to the knapsack set $\text{conv}(\{z\in\{0,1\}^n:\sum_{i=1}^na_iz_i\le b\})$.  If we define $f(x)=g\left(a^\top x\right)$, where $g(t)=\delta\left(t;\left\{t\in\R:t\le b \right\}\right)$, it is clear that $\set Q_B=\R_+\times\set T$, which means $\set T$ is isomorphic to one facet of $\text{conv}(\set Q_B)$. {Since any compact polyhedral description of $\text{conv}(\set Q_B)$ would imply the existence of a polynomial time algorithm to optimize over $\set T$, we conclude that such description is not possible unless $\set P=\np$}.
\subsection*{{Extensions for Rank-one quadratic programs}}
{In this section, we discuss two natural extensions for rank-one quadratic programs. We first consider the extension that box-constraints are added for continuous variables. Specifically, consider} the following mixed-integer quadratic program 
\begin{equation}\label{eq:boxed-QP}
	\begin{aligned}
		\min_{x,z}\;&\left(a^\top x \right)^2+b^\top x+c^\top z\\
		\text{s.t. }&0\le x\le z, \\
		&z\in\{0,1\}^n.
	\end{aligned}
\end{equation}
We aim to show \eqref{eq:boxed-QP} is \np-hard in general. To achieve this goal, we show that \eqref{eq:boxed-QP} includes the following well known 0-1 knapsack problem \eqref{eq:knapsack} as a special case. 
\begin{equation}\label{eq:knapsack}
	\begin{aligned}
		\min_{z\in\{0,1\}^n}\;&-v^\top z\\
		\text{s.t. }&w^\top z\le W,
	\end{aligned}
\end{equation}
where $(v,w,W)\in\Z_+^{2n+1}$ are nonnegative integers such that $w_i\le W\le \sum_j w_j\;\forall i\in[n]$.
\begin{proposition}\label{prop:knapsack}
	The knapsack problem \eqref{eq:knapsack} is equivalent to the optimization problem
	\begin{equation}\label{eq:knapsack-QP}
		\begin{aligned}
			\min_{x,z}\;&M_1\left(Wx_0+\sum_{i=1}^{n}w_ix_i-W\right)^2-M_2\left(\sum_{i=1}^n x_i\right)+\sum_{i=1}^{n}(M_2-v_i)z_i\\
			\text{s.t. }&0\le x_i\le z_i,\; i=0,1,\dots,n\\
			&z_i\in\{0,1\},\;i=0,1,\dots,n,
		\end{aligned}
	\end{equation}where $M_1=\sum_{i=1}^{n}v_i+1$ and $M_2=2nW^2M_1+1$ are polynomial in the input size.
\end{proposition}
\begin{proof}
	Denote the objective function by $\theta(x,z)$. First, since the objective function does not involve $z_0$, $z_0$ can be safely taken as $1$. Second, we now show that $M_2$ is large enough to force $x_i=z_i$ for all $i\in [n]$. Specifically, for any $i\in[n]$, since $x_0,x_i\le 1$ and $w_i\le W$, it holds that
	\begin{align*}
		\partialDerivative{\theta(x,z)}{x_i}&=2M_1w_i\left(Wx_0+\sum_{i=1}^{n}w_ix_i-W\right)-M_2\\
		&\le 2nW^2M_1-M_2<0.
	\end{align*}
	That is, $\theta(x,z)$ is decreasing with respect to $x_i,i\in[n]$, which implies $x_i=z_i$ in any optimal solution. It follows that \eqref{eq:knapsack-QP} can be simplified to
	\begin{equation}\label{eq:knapsack-augmented-QP}
		\begin{aligned}
			\min_{x_0,z}\;&M_1\left(Wx_0+\sum_{i=1}^{n}w_iz_i-W\right)^2-\sum_{i=1}^{n}v_iz_i\\
			\text{s.t. }
			&0\le x_0\le1,\;z_i\in\{0,1\},\;i=1,\dots,n,
		\end{aligned}
	\end{equation}
	
	Next, we claim that $M_1$ is large enough to ensure that the optimal solution to \eqref{eq:knapsack-augmented-QP} satisfies $Wx_0+\sum_{i=1}^{n}w_iz_i-W=0,$ that is, $ x_0=1-\left(w^\top z\right)/W$. To prove it rigorously, observe that the minimum value of \eqref{eq:knapsack-augmented-QP} must be non-positive since $x_0=1,z=0$ is a feasible solution with objective value equal to $0$. Moreover, since $1-\left(w^\top z\right)/W\le 1$, if $x_0\neq1-\left(w^\top z\right)/W$ at the optimal solution to \eqref{eq:knapsack-augmented-QP}, then $1-\left(w^\top z\right)/W<0$ and the optimal $x_0$ must attain its lower bound $0$. Furthermore, $1-\left(w^\top z\right)/W<0$ implies that $w^\top z\ge W+1$ since $w,\;W$ and $z$ are nonnegative integers. In this case, setting $x_0=0$, the minimum objective value of \eqref{eq:knapsack-augmented-QP} can be written as \begin{align*}
		M_1\left(\sum_{i=1}^{n}w_iz_i-W\right)^2-\sum_{i=1}^{n}v_iz_i>M_1-\sum_i v_i=1>0,
	\end{align*} 
	which contradicts the non-positivity of the optimal objective value.

	Therefore, we can substitute out $x_0=1-\left(w^\top z\right)/W$ and \eqref{eq:knapsack-QP} further reduces to 
	\begin{equation*}
		\begin{aligned}
			\min_{z}\;&-\sum_{i=1}^{n}v_iz_i\\
			\text{s.t. }
			&0\le 1-\frac{w^\top z}{W}\le1,\\
			&z_i\in\{0,1\},\;i=1,\dots,n,
		\end{aligned}
	\end{equation*}
	which  is equivalent to \eqref{eq:knapsack}, because $\left(w^\top z\right)/W\ge 0$ holds trivially.
\end{proof}
Due to the equivalence between optimization problems and separation problems \cite{grotschel1981ellipsoid}, Proposition~\ref{prop:knapsack} indicates that it is impossible to extend the analysis in Section~\ref{sec:extended} and Section~\ref{sec:rank-one} to the case with bounded continuous variables. 

{Besides the addition of box-constraints, another natural extension for rank-one quadratic programs is the ``rank-one+diagonal" objective. % Next, we turn to  extension for rank-one quadratic programs, where the quadratic form in the objective is assumed to be ``rank-one+diagonal".
	 Precisely, consider
\begin{equation}\label{eq:diag+rank-one}
	\begin{aligned}
		\min_{x,z}\;&\left(a^\top x \right)^2+\sum_{i=1}^nd_ix_i^2+b^\top x+c^\top z\\
		\text{s.t. }&x\circ(1-z)=0, \\
		&x\in\R^n_+,\;z\in\{0,1\}^n,
	\end{aligned}
\end{equation}
where $a,b,c\in\R^n$, $d\in\R^n_+$, and vector $u\circ v$ represents the Hadamard product of $u$ and $v$, i.e. $(u\circ v)_i=u_iv_i$. \citet{HGA:2x2} proves that problem \eqref{eq:diag+rank-one} includes the \np-hard subset sum problem as a special case.
\begin{proposition}[\cite{HGA:2x2}]
	Problem~\eqref{eq:diag+rank-one} is \np-hard.
\end{proposition}
For the same reason mentioned above, it is unlikely to get a tractable convex hull description of the corresponding mixed-integer epigraph \[\left\{(t,x,z)\in\R\times\R_+^n\times\{0,1\}^n:\;t\ge \left(a^\top x \right)^2+\sum_{i=1}^nd_ix_i^2,\;x\circ(1-z)=0 \right\}.\] 
}

{
\section{Experiments with signal denoising problems}\label{sec:experiments}
In this section, we test the effectiveness of the low-rank formulations derived in this paper on signal denoising problems. In particular, given noisy observations $c\in \R^n$ of a temporal process $X\in \R^n$, the goal is to recover the true values of the process. We assume (i) that the temporal process is sparse, i.e., $\|X\|_0$ is small; (ii) most observations are noisy realizations of $X$, i.e., $(X_i-c_i)^2$ is small; (iii) a small (and unknown) subset $J$ of observations have been subjected to arbitrarily large corruptions, thus the values $c_j$ for $j\in J$ can be arbitrarily far from $X_j$. The problem of inferring the true values of the process $X$ can be mathematically formulated as 
\begin{equation}\label{eq:signalDenoising}
	\begin{aligned}
		\min_{x,v,z,w}\;&\sum_{i=1}^n(x_i-v_i-c_i)^2+\Omega\sum_{i=\ell+1}^n\left( x_i-\sum_{j=1}^\ell \alpha^j x_{i-\ell+j-1} \right)^2\\
		\text{s.t. }&x\circ(1-z)=0,\; \sum_{i=1}^n z_i\le k_1,\; z\in\{0,1\}^n\\
		&v\circ(1-w)=0,\;\sum_{i=1}^nw_i\le k_2,\;w\in\{0,1\}^n,
	\end{aligned}
\end{equation} 
where $k_1$ is the estimated number of nonzero elements in the true signal, $k_2$ is the maximum number of outliers in data $c$, $\Omega>0$ is the weight of the smoothness part in the objective, $0<\alpha<1$, and $\ell$ is the length of the convolution kernel.

Note that the objective function consists of two parts -- the first sum encodes the fitness of the estimator to the observations; the second sum incorporates the smoothness prior of the true signal that the current signal value $x_i$ is assumed to be similar to those incurred before the time stamp $i$ and the proximity decays in temporal distance at rate $\alpha$. Observe that $w_i$ indicates whether the observation $i$ is an outlier. Indeed, if $w_i=0$, then $v_i=0$ by the complementarity constraint $v_i(1-w_i)=0$ and $c_i$ is used in \eqref{eq:signalDenoising}. On the other hand, if $w_i=1$, then at the optimal solution one must have $v_i=x_i-c_i$. In this case, $c_i$ plays no role in the objective and datapoint $i$ is thus discarded as an outlier.  

We point out that problem \eqref{eq:signalDenoising} is closely related to several problems in the literature. On one hand, if $\ell=1$ and $\alpha=1$, then the smoothness terms reduce to $(x_i-x_{i-1})^2$. In this case, problems with sparse signals and no outliers (i.e., $k_2=n$) where considered in \cite{atamturk2021sparse}, and problems with outliers but no sparsity (i.e., $k_1=n$) were considered in \cite{gomez2021outlier}. Moreover, if $\ell=1$, then the penalized version of \eqref{eq:signalDenoising}, where cardinality constraints on sparsity and outliers are replaced with fixed costs, is polynomial-time solvable \cite{han2022polynomial}. On the other hand, if the width $\ell$ is large and terms $\alpha^j$ are replaced with arbitrary constants, then \eqref{eq:signalDenoising} is closely related with sparse/robust versions of linear regression. In particular, there have been several recent papers related to sparse linear regression \cite{atamturk2020safe,bertsimas2016best,cozad2014learning,hazimeh2022sparse,xie2020scalable}, but relatively few mixed-integer optimization approaches for problems with outliers \cite{zioutas2005deleting,zioutas2009quadratic} or both outliers and sparsity \cite{insolia2022simultaneous}.

\subsection{Instance generation}
In this section, we describe how we generate test instances. Set $\alpha = 0.9$, $k_1=3n/50$ and $k_2=n/100$. The noised signal $c$ is produced in the following way:
\begin{enumerate}
	\item Initially, $c=0$
	\item Repeat the following process to generate $s=n/50$ ``spikes" of the signal with each spike consisting of $h=5$ non-zero values 
	\begin{enumerate}
		\item Select an index $i_0$ uniformly between 1 and $n+1-h$ as the starting position of a given spike.
		\item Draw an vector $v\in\R^h$ from a multivariate normal distribution $\mathcal{N}(0,\Sigma)$, where $\Sigma$ is the covariate matrix with $\Sigma_{ij}=\frac{i(h+1-j)}{h+1}$.
		\item Update $c_{i_0+i}\leftarrow c_{i_0+i}+v_i,\; i=1,\dots,h$.
	\end{enumerate}
	\item Add a Gaussian noise  $c\leftarrow c+0.2\varepsilon$, where $\varepsilon\in\R^n$ is sample from the standard multivariate normal distribution $\set N(0,I)$. Then normalize the data $c\leftarrow c/\norm{\infty}{c}$.
	\item Corrupt the data $c$ by repeating the following process $k_2$ times
	\begin{enumerate}
		\item Select an index $i$ uniformly between 1 and $n$.
		\item Update $c_i\leftarrow c_i+4\epsilon$, where $\epsilon$ is drawn from the standard normal distribution.
	\end{enumerate}
\end{enumerate}
All experiments in this section are conducted using Gurobi~9.0 solver with default settings on a laptop with a 2.30GHz $\text{Intel}^\text{\textregistered}$
$\text{Core}^{\text{\tiny TM}}$ i9-9880H CPU and 64 GB main memory. The time limit is set as 600 seconds.

\subsection{Formulations}
Three mixed-integer reformulations of \eqref{eq:signalDenoising} are implemented and compared - $\basic$, $\rankOne$, and $\rankTwo$.

\subsubsection{$\basic$} In $\basic$, we reformulate the complementarity constraints using big-M formulations where $M=10^4$
\begin{equation}\tag{$\basic$}
	\begin{aligned}
		\min_{x,v,z,w}\;&\sum_{i=1}^n(x_i-v_i-c_i)^2+\Omega\sum_{i=\ell+1}^{n}\left( x_i-\sum_{j=1}^{\ell}\alpha^{j}x_{i-\ell+j-1} \right)^2\\
		\text{s.t. }&-Mz\le x\le Mz,\;\sum_i z_i\le k_1, \;z\in\{0,1\}^n\\
		&-Mw\le v\le Mw,\;\sum_i w_i\le k_2,\; w\in\{0,1\}^n.
	\end{aligned}	
\end{equation}
We comment that the continuous relaxation of $\basic$ --the optimization problem obtained by dropping the integrality constraints of $x$ and $w$-- always yields a trivial optimal value $0$. Indeed, since $M$ is sufficiently large, one can take $v=c$, $x=z=0$ and $w=v/M$ as a feasible solution to the relaxation of $\basic$ whose objective is simply $0$.  
\subsubsection{$\rankOne$} In $\rankOne$, each individual fitness term and smoothness term in \eqref{eq:signalDenoising} are replaced with the conic representation of the rank-one formulation suggested by  Proposition~\ref{prop:ran-one-explicit}
\begin{equation}\tag{\rankOne}
	\begin{aligned}
		\min_{x,v,z,w}\;&\sum_{i=1}^{n}t_i+\Omega\sum_{i=\ell+1}^{n}s_i-2\sum_{i=1}^nc_i(x_i-v_i)+\norm{2}{c}^2\\
		\text{s.t. }& t_i\ge (x_i-v_i)^2,\; t_i(z_i+w_i)\ge (x_i-v_i)^2&i=1,\dots,n\\
		&s_i\ge\left( x_i-\sum_{j=1}^{\ell}\alpha^{j}x_{i-\ell+j-1} \right)^2&i=\ell+1,\dots,n\\
		&s_i\left(\sum_{j=i-\ell}^{i}z_j\right)\ge\left( x_i-\sum_{j=1}^{\ell}\alpha^{j}x_{i-\ell+j-1} \right)^2& i=\ell+1,\dots,n \\
		&-Mz\le x\le Mz,\;\sum_i z_i\le k_1, \;z\in\{0,1\}^n\\
		&-Mw\le v\le Mw,\;\sum_i w_i\le k_2,\; w\in\{0,1\}^n.
	\end{aligned}	
\end{equation}
\subsubsection{$\rankTwo$} In $\rankTwo$, for each $i=1,\dots,\ell$, we replace the $i^{\text{th}}$ fitness term with its rank-one strengthening. For each $i=\ell+1,\dots,n$, we combine the $i^{\text{th}}$ fitness and smoothness term and replace their sum with the rank-two strengthening. Specifically, define the rank-two mixed-integer set \[\set Q_2=\left\{(t,x,z,v,w,y,u):\begin{aligned}
	&t\ge(y-v)^2+\Omega\left(y-\sum_{i=1}^\ell x_i\right)^2\\
	&x\circ(1-z)=0,\;v(1-w)=0,\;y(1-u)=0\\
	&z\in\{0,1\}^\ell,w\in\{0,1\},u\in\{0,1\}
\end{aligned}\right\}.\]
The description of $\conv\set (Q_2)$ in the lifted space can be readily derived from Theorem~\ref{thm:main} using disjunctive programming techniques, which we give in Proposition~\ref{prop:signal} in the appendix. Then $\rankTwo$ is given as
\begin{equation}\tag{$\rankTwo$}
	\begin{aligned}
		\min_{x,v,z,w}&\sum_{i=1}^nt_i-2\sum_{i=1}^nc_i(x_i-v_i)+\norm{2}{c}^2\\
		\text{s.t. }& t_i\ge (x_i-v_i)^2, t_i(z_i+w_i)\ge (x_i-v_i)^2&i=1,\dots,\ell\\
		&(t_i,a\circ x_{(i-\ell):(\ell-1)},z_{(i-\ell):(\ell-1)},v_i,w_i,x_i,z_i)\in\conv(\set Q_2)&i=\ell+1,\dots,n\\
		&-Mz\le x\le Mz,\;\sum_i z_i\le k_1, \;z\in\{0,1\}^n\\
		&-Mw\le v\le Mw,\;\sum_i w_i\le k_2,\; w\in\{0,1\}^n,
	\end{aligned}
\end{equation}
where $x_{i:j}$ represents the subvector of $x$ indexed from $i$ to $j$ and $a=(\alpha^\ell,\alpha^{\ell-1},\dots,\alpha)$. 
\subsection{Results}
Each row of Table~\ref{tab:kernel-width} and Table~\ref{tab:weight} shows an average statistics over five instances with the same parameters. The tables display the dimension of the problem $n$, the initial gap (IGap), the solution time for solving the corresponding mixed-integer formulation (Time), the end gap provided by the solver at termination (EGap), the number of nodes explored by the solver (Nodes), the number of instances solved to optimality within the time limit (\#), and the root improvement of $\rankTwo$ compared with $\basic$ (RI-basic) and $\rankOne$ (RI-rankOne).  The initial gap is computed as $\text{IGap}=\frac{\text{obj}_{\text{best}}-\text{obj}_{\text{cont}}}{|\text{obj}_{\text{best}}|}\times 100$, where $\text{obj}_{\text{best}}$ is the objective value of the feasible solution found and $\text{obj}_{\text{cont}}$ represents the objective value of the continuous relaxation, that is, obtained by relaxing $\{0,1\}^n$ to $[0,1]^n$ in the three formulations. The root improvement is computed as $\text{RI-basic}=\frac{\text{obj}_{\text{rankTwo}}-\text{obj}_{\text{basic}}}{\text{obj}_{\text{best}}-\text{obj}_{\text{basic}}}\times100$, where $\text{obj}_{\text{rankTwo}}$ is the objective value of the convex relaxation of $\rankTwo$ and  where $\text{obj}_{\text{basic}}$ is the objective value of the convex relaxation of $\basic$. The computation of RI-rankOne is similar. Finally, we remark that the initial gap of $\basic$ is not reported because it is always nearly 0, which is expected and explained above.

Table~\ref{tab:kernel-width}  displays the computational results for varying $\ell\in\{1,2,5\}$ and fixed $\Omega=0.05$. As we expect, $\basic$ has the worst performance and can only solve three instances in  the low-dimensional setting $n=100$. On the other hand, $\rankTwo$ is able to solve almost all instances to optimality when $\ell=1,2$ by exploring much less nodes in the branch-and-bound  tree regardless of the problem dimension. The root improvement of $\rankTwo$
over $\basic$ is significant in all settings. The performance of $\rankOne$ is in the middle of $\basic$ and $\rankTwo$ - it can solve all instances in low-dimension settings and none when $n=200,300$.  In addition, RI-rankOne ranges from $8\%$ to $50\%$ and achieves the best $50\%$ as  $\ell=2$ and $n=100$. Table~\ref{tab:weight} displays the computational results for fixed $\ell=2$ and varying $\Omega\in\{0.01,0.05,0.25,0.5\}$. The conclusions we can draw from Table~\ref{tab:weight} are similar to the ones from Table~\ref{tab:kernel-width}. It is worthing noting that in the case where a certain instance can be solved by both $\rankOne$ and $\rankTwo$, the solution time of $\rankOne$ outperforms $\rankTwo$. However, many instances that are not solvable using $\rankOne$ can be solved by $\rankTwo$ especially in more difficult settings. In summary, we conclude that $\rankTwo$ is preferable to utilizing for challenging instances of the signal denoising problem \eqref{eq:signalDenoising}.

The computational profile of three formulations on all 105 instances is presented in Figure~\ref{fig:profile_signal}. Among them, $\rankOne$ successfully solves 41 instances in 600 seconds. In comparison,  $\rankTwo$ is able to solve the same number of instances in 5.7 seconds, which is two orders of magnitude faster than $\rankOne$. Furthermore, $\rankTwo$ surpasses $\rankOne$ by solving 81 instances within the given time limit, which is approximately twice the number of instances solved by $\rankOne$. Overall, these findings demonstrate the superior performance of $\rankTwo$ relative to $\rankOne$ in terms of solvability. 
}

\section{Conclusions}\label{sec:conclusions}
In this paper, we propose a new disjunctive programming representation of the convex envelope of a low-rank convex function with indicator variables and complementary constraints. The ensuing formulations are substantially more compact than alternative disjunctive programming formulations. As a result, it is substantially easy to project out the additional variables to recover formulations in the original space of variables, and to implement the formulations using off-the-shelf solvers. { Moreover, the computational results evidently demonstrate the efficacy of the proposed formulations when applied to well-structured problems, particularly in problems where the low rank functions are involving few variables as well. }

{
\begin{figure}[h!]
	\includegraphics[trim={11cm 6cm 10cm 6cm}, clip,width=\textwidth]{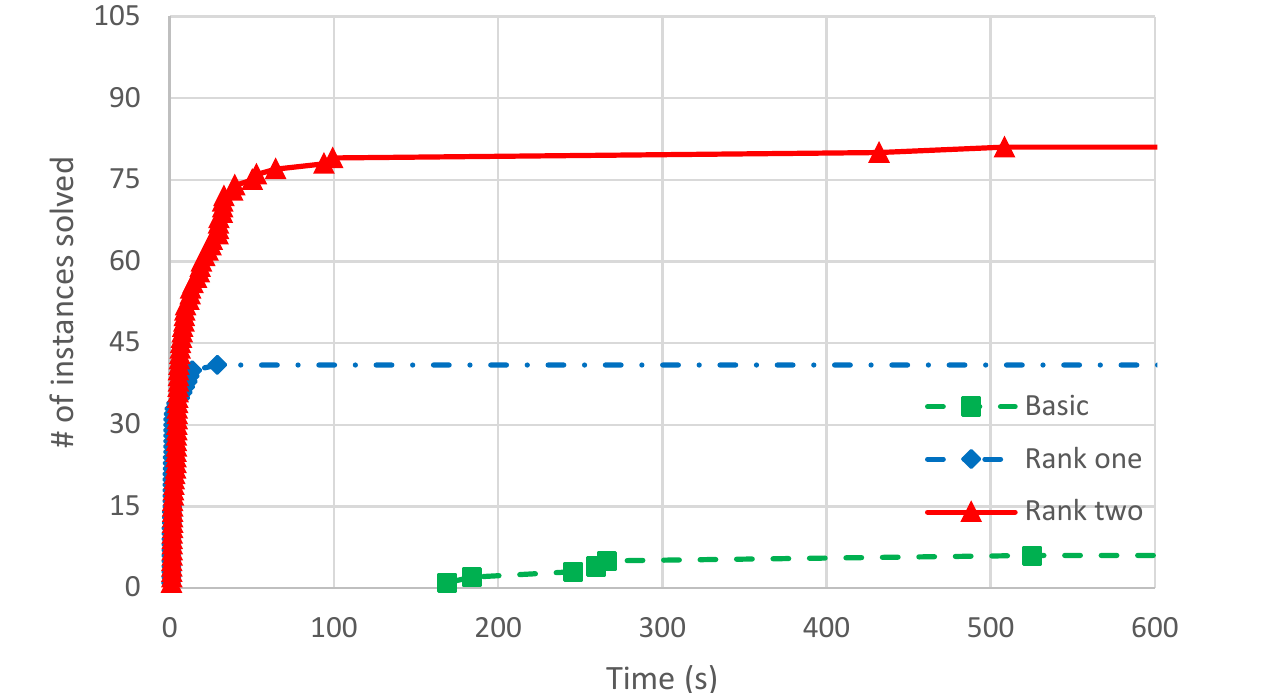}
	\caption{Number of instances solved as a function of time}\label{fig:profile_signal}
\end{figure}
\begin{table}[h]
	%\footnotesize
	\centering
	%\small
	\caption{Varying kernel width $\ell$ and fixed $\Omega=0.05$}\label{tab:kernel-width}
	\resizebox{\linewidth}{!}{
		\begin{tabular}{rr|rrrc|rrrrc|rrrrrrc}
			\toprule		
			\multirow{2}{*}{$n$} & \multirow{2}{*}{$\ell$}&\multicolumn{4}{c|}{\basic}&\multicolumn{5}{c|}{$\rankOne$}&\multicolumn{7}{c}{$\rankTwo$} \\\cmidrule(lr){3-6} \cmidrule(lr){7-11} \cmidrule(lr){12-18}
			& & EGap & Time & Nodes & \# &   IGap &   EGap &   Time &  Nodes & \# &   IGap &   EGap &   Time &  Nodes & \# & RI-basic &   RI-rankOne\\
			\midrule
			100 &          1 &      49.00 &        447 &    180,297 &          2 &       4.30 &       0.00 &       1.13 &        342 &          5 &       2.09 &       0.00 &       1.65 &        7.8 &          5 &      97.91 &      50.37 \\
			
			100 &          2 &      54.40 &        585 &    182,314 &          1 &       5.68 &       0.00 &       1.36 &        131 &          5 &       4.22 &       0.00 &       1.51 &        3.8 &          5 &      95.78 &      26.02 \\
			
			100 &          5 &      61.01 &        600 &     39,610 &          0 &       8.99 &       0.00 &       2.91 &        344 &          5 &       7.21 &       0.00 &       5.60 &       16.8 &          5 &      92.79 &      20.42 \\\midrule
			
			200 &          1 &      43.30 &        600 &      3,394 &          0 &       3.91 &       3.03 &        600 &      9,279 &          0 &       2.07 &       0.00 &       23.5 &        155 &          5 &      97.93 &      48.08 \\
			
			200 &          2 &      51.22 &        600 &      3,391 &          0 &      10.36 &       9.45 &        600 &      6,610 &          0 &       8.18 &       0.00 &       32.3 &      102.8 &          5 &      91.82 &      21.31 \\
			
			200 &          5 &      56.95 &        600 &      3,688 &          0 &      10.82 &      10.31 &        600 &      7,986 &          0 &       9.70 &      10.91 &        600 &        531 &          0 &      90.30 &      11.13 \\\midrule
			
			300 &          1 &      82.17 &        600 &      3,149 &          0 &       4.78 &       7.94 &        600 &      3,414 &          0 &       2.50 &       0.00 &       23.4 &         58 &          5 &      97.50 &      48.82 \\
			
			300 &          2 &      87.06 &        600 &      4,288 &          0 &       8.70 &      10.09 &        600 &      2,701 &          0 &       7.18 &       4.02 &        141 &        311 &          4 &      92.82 &      17.99 \\
			
			300 &          5 &      94.55 &        600 &      2,520 &          0 &      17.33 &      18.31 &        600 &      2,766 &          0 &      16.00 &      31.18 &        600 &      292.4 &          0 &      84.00 &       8.68 \\
			\bottomrule
		\end{tabular}  
	}
\end{table}

\begin{table}[h]
	%\footnotesize
	\centering
	%\small
	\caption{Varying $\Omega$ and fixed kernel-width $\ell=2$}\label{tab:weight}
	\resizebox{\linewidth}{!}{
		\begin{tabular}{rr|rrrc|rrrrc|rrrrrrc}
			\toprule		
			\multirow{2}{*}{$n$} & \multirow{2}{*}{$\Omega$}&\multicolumn{4}{c|}{$\basic$}&\multicolumn{5}{c|}{$\rankOne$}&\multicolumn{7}{c}{$\rankTwo$} \\\cmidrule(lr){3-6} \cmidrule(lr){7-11} \cmidrule(lr){12-18}
			& & EGap & Time & Nodes & \# &   IGap &   EGap &   Time &  Nodes & \# &   IGap &   EGap &   Time &  Nodes & \# & RI-basic &   RI-rankOne\\
			\midrule	
			100 &       0.01 &      94.42 &        600 &    146,503 &          0 &       1.38 &       0.00 &       1.32 &        159 &          5 &       0.90 &       0.00 &       1.67 &          3 &          5 &      99.10 &      35.50 \\
			
			100 &       0.05 &      26.61 &        446 &     50,971 &          2 &       7.31 &       0.00 &       1.22 &        253 &          5 &       4.92 &       0.00 &       1.95 &          8 &          5 &      95.08 &      32.67 \\
			
			100 &       0.25 &      41.19 &        600 &     65,345 &          0 &      17.50 &       0.00 &       1.76 &        740 &          5 &      13.33 &       0.00 &       5.41 &         36 &          5 &      86.67 &      25.46 \\
			
			100 &       0.50 &      18.18 &        532 &     70,338 &          1 &      28.99 &       0.00 &       2.96 &      1,056 &          5 &      24.14 &       0.00 &       4.21 &         23 &          5 &      75.86 &      17.01 \\\midrule
			
			200 &       0.01 &      34.31 &        600 &      4,653 &          0 &       1.92 &       0.00 &       9.98 &      1,438 &          5 &       1.43 &       0.00 &       5.81 &         11 &          5 &      98.57 &      26.23 \\
			
			200 &       0.05 &      42.16 &        600 &      4,811 &          0 &       7.01 &       5.55 &        600 &      7,430 &          0 &       5.27 &       0.75 &        134 &        244 &          4 &      94.73 &      25.56 \\
			
			200 &       0.25 &      61.65 &        600 &      3,262 &          0 &      20.11 &      20.25 &        600 &      6,573 &          0 &      15.60 &       3.72 &        217 &        579 &          4 &      84.40 &      23.02 \\
			
			200 &       0.50 &      75.17 &        600 &      3,414 &          0 &      30.80 &      29.86 &        600 &      6,092 &          0 &      24.38 &      17.77 &        484 &        614 &          1 &      75.62 &      22.39 \\\midrule
			
			300 &       0.01 &      78.96 &        600 &      2,785 &          0 &       1.65 &       1.76 &        486 &      4,120 &          1 &       1.15 &       0.00 &       29.3 &         51 &          5 &      98.85 &      30.19 \\
			
			300 &       0.05 &      87.16 &        600 &      2,959 &          0 &       7.85 &      10.74 &        600 &      2,844 &          0 &       5.77 &       0.00 &        122 &        339 &          5 &      94.23 &      26.58 \\
			
			300 &       0.25 &      80.01 &        600 &      3,192 &          0 &      20.46 &      20.87 &        600 &      2,780 &          0 &      16.12 &      15.93 &        368 &        369 &          2 &      83.88 &      21.79 \\
			
			300 &       0.50 &      87.68 &        600 &      3,265 &          0 &      31.09 &      31.80 &        600 &      2,485 &          0 &      25.20 &      25.51 &        488 &        543 &          1 &      74.80 &      19.17 \\
			\bottomrule
		\end{tabular}  
	}
\end{table}
}

\section*{Acknowledgments}

This research is supported in part by the National Science Foundation under grant CIF 2006762.
\linespread{1.0}

\bibliographystyle{apalike}
\bibliography{Bibliography}

\appendix

\section{On computational experiments with branch-and-bound}

As mentioned in \S\ref{sec:implementation}, the extended formulation \eqref{eq:conicQuadraticRepresentability} did not produce good results when used in conjunction with CPLEX branch-and-bound solver. To illustrate this phenomenon, Table~\ref{tab:branch-and-bound} shows details on the performance of the solver in a single representative instance, but similar behavior was observed in \emph{all} instances tested. The table shows from left to right: the time required to solve the convex relaxation via interior point methods and the lower bound produced by this relaxation (note that this is not part of the branch-and-bound algorithm); the time required to process the root node of the branch-and-bound tree, and the corresponding lower bound obtained; the time used to process the branch-and-bound tree, the number of branch-and-bound nodes explored, and the lower bound found after processing the tree (we set a time limit of 10 minutes). 

\begin{table}[!h]
	\caption{Performance of CPLEX solver in an instance with $n=500$ and $k=10$. Default settings are used, and a time limit of 10 minutes is set. The optimal objective value in the particular instance is 1.47.}
	\label{tab:branch-and-bound}
	\begin{tabular}{c |c c | c c | c c c}
		\hline
		\multirow{2}{*}{Method}& \multicolumn{2}{c|}{\underline{Convex relaxation}}&\multicolumn{2}{c|}{\underline{Root node}}&\multicolumn{3}{c}{\underline{Branch-and-bound}}\\
		&Time(s)&LB&Time(s)&LB&Time(s)&Nodes&LB\\
		\hline
		Without \eqref{eq:conicQuadraticRepresentability} & 0.2 & 1.09 &3.5 & 1.20 &7.7&460&1.47\\
		With \eqref{eq:conicQuadraticRepresentability} & 2.4 & 1.30 & 45.9 & 0.90 &600 & 4,073 & 0.99\\
		\hline
	\end{tabular}
\end{table}

Our expectations a priori were as follows: using inequalities \eqref{eq:conicQuadraticRepresentability} should result in harder convex relaxations solved, and thus less nodes explored within a given time limit; on the other hand, due to improved relaxations and higher-quality lower bounds, the algorithm should be able to prove optimality after exploring substantially less nodes. Thus, there should be a tradeoff between the number of nodes to be explored and the time required to process each node. From Table~\ref{tab:branch-and-bound}, we see that there is no tradeoff in practice. 

The performance of the solver without inequalities \eqref{eq:conicQuadraticRepresentability} is as expected. While just solving the convex relaxation via interior point methods requires 0.2 seconds, there is an overhead of 3 seconds to process the root node due to preprocessing/cuts/heuristics and additional methods used by the solver (and the quality of the lower bound at the root node is slightly improved as a result).Then, after an additional 4 seconds used to explore 460 nodes, optimality is proven. 

The performance of the solver using inequalities \eqref{eq:conicQuadraticRepresentability} defied our expectations. In theory, the more difficult convex relaxation can be solved with an overhead of 2 seconds, resulting in a root improvement of $(1.30-1.09)/(1.47-1.09)=55.3\%$ (better than the one achieved by default CPLEX). In practice, the overhead is 40 seconds, and results in a \emph{degradation} of the relaxation, that is, the lower bound proved at the root node is worse than the natural convex relaxation of problem without inequalities \eqref{eq:conicQuadraticRepresentability}. From that point out, the branch-and-bound progresses slowly due to the more difficult relaxations, and the lower bounds are worse throughout the tree. Even after the time limit of 10 minutes and over 4,000 nodes explored, the lower bound proved by the algorithm is still worse than the natural convex relaxation. 

While we cannot be sure about the exact reason of this behavior, we now make an educated guess. Most conic quadratic branch-and-bound solvers such as CPLEX do not use interior point methods to solve relaxations at each node of the branch-and-bound tree, but rather rely on polyhedral outer approximations in an extended space to benefit from the warm-start capabilities of the simplex method. We conjecture that while formulation \eqref{eq:conicQuadraticRepresentability} might not be particularly challenging to solve via interior point methods, it might be difficult to construct a good-quality outer approximation of reasonable size. If so, then the actual relaxation used by the solver is possibly a poor-quality linear outer approximation of the feasible region induced by \eqref{eq:conicQuadraticRepresentability}, and is still difficult to solve, resulting in the worst of both worlds. We point out that the second author has encountered similar counterintuitive behavior with solvers (other than CPLEX) based on linear outer approximations \cite{atamturk2018strong}. {We also attempted that setting the MIQCP strategy switch (parameter MIQCPStrat) to one, which ``tells CPLEXs to solve a QCP relaxation of the model at each node". However, in such cases, the optimization terminates with numerical issues in most cases.}

{
\section{Convex hull description of $\set Q_2$}
In this section, we give the convex hull description of $\set Q_2$ in the lifted space. 
We first prove a lemma to simplify the notation. Consider  
\[ \set S_i=\left\{(t,x,z):t\ge f_i(x_{\set I_i}),\;z_j=1\;\forall j\in\set I_i, \;x_j=0,z_j\in\{0,1\}\;\forall j\not\in\set I_i  \right\}\;i\in[n], \]
where each $\set I_i\subseteq[n]$ and $f_i$ is a closed convex function with $f_i(0)=0$. 
\begin{lemma}\label{lem:auxilliary-elimination}
	Let $\set S=\bigcup_i \set S_i$. A point $(t,x,z)\in\conv\set S$ if and only if the following system is consistent
	\begin{align*}
		&t=\sum_i t_i\\
		&t_i\ge f^\pi(x^i_{\set I_i};\lambda_i)&\forall i\in[n]\\
		&x_i=\sum_{j:i\in \set I_j} x^j_i&\forall i\in[n]\\
		&\sum_{j:i\in\set I_j}\lambda_ j\le z_i \le 1&\forall i\in[n]\\
		&\lambda\ge0,\sum_{i\in[n]}\lambda_i=1.
	\end{align*}
\end{lemma}
\begin{proof}
	Using disjunctive programming techniques, one can see that $(t,x,z)\in\conv\set S$ if and only if the following system is consistent
	\begin{align*}
		&t=\sum_i t_i\\
		&x_i=\sum_j x^j_i,\;z_i=\sum_jz^j_i&\forall i\in[n]\\ 
		&x_j^i=0, z_j^i\in[0,\lambda_i]\;&\forall j\notin\set I_i,\forall i\in[n]\\
		&z_j^i=\lambda_i\;&\forall j\in\set I_i,\forall i\in[n]\\
		&t_i\ge f^\pi(x^i_{\set I_i};\lambda_i)&\forall i\in[n]\\
		&\lambda\ge0,\sum_{i\in[n]}\lambda_i=1.
	\end{align*}
First observe that one can easily substitute out $x^i_j\;\forall j\notin\set I_i$ and $z^i_j\;\forall j\in\set I_i$ in above expression. What is more, since for each $j\notin\set I_i$, $z_j^i\in[0,\lambda_i]$ is independent of all other variables except $z_i$ and $\lambda_i$. Thus, one can use Fourier-Motzkin eliminating to substitute such $z_j^i$ with their bounds, which gives the upper bound of $z_i\le \sum_i\lambda_i=1$.
\end{proof}
Recall the definition of \[\set Q_2=\left\{(t,x,z,v,w,y,u):\begin{aligned}
	&t\ge(y-v)^2+\Omega\left(y-\sum_{i=1}^\ell x_i\right)^2\\
	&x\circ(1-z)=0,\;v(1-w)=0,\;y(1-u)=0\\
	&x\in\R^\ell, z\in\{0,1\}^\ell,w\in\{0,1\},u\in\{0,1\}
\end{aligned}\right\}.\]
\begin{proposition}\label{prop:signal}
	The closed convex hull of $\set Q_2$ is given by the following extended formulation
	\begin{align*}
		&t\ge \sum_{i}( t^{\set V_i}+t^{\set V_{iy}}+t^{\set V_{iv}})+t^{\set V_{y}}+t^{\set V_v}+t^{\set V_{yv}}\\
		&x_i = x^{\set V_i}_i+x_{i}^{\set V_{iy}}+x_{i}^{\set V_{iv}}+x_i^{\set R}&i\in[\ell]\\
		&y = y^{\set V_y}+\sum_{i=1}^\ell y^{\set V_{iy}}+y^{\set V_{yv}}+y^{\set R}\\
		&v=v^{\set V_v}+\sum_{i=1}^\ell v^{\set V_{iv}}+v^{\set V_{yv}}+v^{\set R}\\
		&y^{\set R}=v^{\set R}=\sum_{i=1}^\ell x_i^{\set R}\\
		&\lambda^{\set V_i}+\lambda^{\set V_{iy}}+\lambda^{\set V_{iv}}+\lambda^{\set R}\le z_i\le 1&\forall i\in[\ell]\\
		&\lambda^{\set V_y}+\sum_{i=1}^\ell \lambda^{\set V_{iy}}+\lambda^{\set V_{yv}}+\lambda^{\set R}\le u\le 1\\
		&\lambda^{\set V_v}+\sum_{i=1}^\ell \lambda^{\set V_{iv}}+\lambda^{\set V_{yv}}+\lambda^{\set R}\le w\le 1\\
		&\sum_i\lambda^{\set V_i}+\sum_i\lambda^{\set V_{iy}}+\sum_i\lambda^{\set V_{iv}}+\lambda^{\set V_y}+\lambda^{\set V_v}+\lambda^{\set V_{yv}}+\lambda^{\set V_R}\le 1,\;\lambda\ge0\\
		&t^{\set V_{i}}\lambda^{\set V_i}\ge \Omega \left(x^{\set V_i}\right)^2&\forall i\in[\ell]\\
		&t^{\set V_y}\lambda^{\set V_y}\ge (1+\Omega)\left( y^{\set V_y} \right)^2\\
		&t^{\set V_v}\lambda^{\set V_v}\ge \left( v^{\set V_v} \right)^2\\
		&t^{\set V_{iy}}\lambda^{\set V_{iy}}\ge \left(y^{\set V_{iy}}\right)^2+\Omega\left( y^{\set V_{iy}} -x_i^{\set V_{iy}} \right)^2&\forall i\in[\ell]\\
		&t^{\set V_{iv}}\lambda^{\set V_{iv}}\ge \left(v^{\set V_{iv}}\right)^2+\Omega \left(x_{i}^{\set V_{iv}}\right)^2&\forall i\in[\ell ]\\
		&t^{\set V_{yv}}\lambda^{\set V_{yv}}\ge(y^{\set V_{yv}}-v^{\set V_{yv}})^2+\Omega \left(y^{\set V_{yv}}\right)^2
	\end{align*}
\end{proposition}
\begin{proof}
	Define \begin{align*}
		&\set R=\left\{ (t,x,z,v,w,y,u):t\ge0, y-v=0, y=\sum_jx_j, z_i=w=u=1\;\forall j \right\}\\
		&\set V_0=\left\{ (t,x,z,v,w,y,u):\begin{aligned}
			&t\ge 0, x_j=0, z_j\in\{0,1\}\;\forall j,\\
			&v=y=0, w,u\in\{0,1\}
		\end{aligned}\right\}\\
		&\set V_i=\left\{ (t,x,z,v,w,y,u):\begin{aligned}
			&t\ge \Omega x_i^2,\;z_i=1,\\
			&v=y=0, w,u\in\{0,1\}, x_j=0,z_j\in\{0,1\}\;\forall j\neq i
		\end{aligned}\right\}&\forall i\in[\ell]\\
		&\set V_y=\left\{ (t,x,z,v,w,y,u):\begin{aligned}
			&t\ge (1+\Omega)y^2,u=1, \\
			&v=0,w\in\{0,1\}, x_j=0, z_j\in\{0,1\}\;\forall j
		\end{aligned} \right\}\\
		& \set V_v=\left\{(t,x,z,v,w,y,u):t\ge v^2,w=1, y=x_j=0,u,z_j\in\{0,1\}\;\forall j  \right\}\\
		&\set V_{yv}=\left\{ (t,x,z,v,w,y,u):\begin{aligned}
			&t\ge y^2+\Omega(y-x_i)^2, z_i=u=1, \\
			&v=0,w\in\{0,1\},x_j=0,z_j\in\{0,1\}\;\forall j\neq i
		\end{aligned} \right\}\\
		&\set V_{iy}=\left\{ (t,x,z,v,w,y,u):\begin{aligned}
			&t\ge y^2+\Omega(y-x_i)^2,z_i=u=1,\\
			&v=0,w\in\{0,1\},x_j=0,z_j\in\{0,1\}\;\forall j\neq i 
		\end{aligned}\right\}&\forall i\in[\ell]\\
	&\set V_{iv}=\left\{(t,x,z,v,w,y,u):\begin{aligned}
		&t\ge v^2+\Omega x_i^2, w=z_i=1,\\
		& y=0, u\in\{0,1\}, x_j=0, z_j\in\{0,1\}\;\forall j\neq i
	\end{aligned} \right\}.&\forall i\in[\ell]
	\end{align*}
For a general rank-two set, Theorem~\ref{thm:main} suggest using $\bigO{n^2}$ of pieces to describe its closed convex hull. However, for this special set $\set Q_2$, carefully checking the proof of Theorem~\ref{thm:main} we find that that the above defined pieces are sufficient, i.e. 
\[ \conv\set Q_2=\conv\left( \set V_0\bigcup_i \set V_i\bigcup \set V_y\bigcup\set V_v\bigcup\set V_{yv}\bigcup_i \set V_{iy}\bigcup_i\set V_{iv}\bigcup \set R \right). \] 
The conclusion follows by applying Lemma~\ref{lem:auxilliary-elimination} to the above disjunction and eliminating the variables associated with $\set V_0$. Note that we use superscripts to indicates the piece associated with a certain additional variable in the extended formula of $\conv\set Q_2$.
\end{proof}
Finally, we note that the number of additional variables used in Proposition~\ref{prop:signal} is of order $\bigO{\ell}$.
}

\end{document}

%% file: isomath.tex
\usepackage{amsmath,amsthm,amssymb,enumitem,xspace}

\newcommand{\Z}{\mathbb{Z}}

\newcommand{\R}{\mathbb{R}}

\newcommand{\numberthis}{\addtocounter{equation}{1}\tag{\theequation}}

\DeclareMathOperator{\conv}{conv}

\DeclareMathOperator*{\dom}{dom}

\newcommand{\norm}[2]{\left\|#2\right\|_{#1}}

\newcommand{\partialDerivative}[2]{\frac{\partial#1}{\partial #2}}

\newcommand{\bigO}[1]{\ensuremath{\mathcal{O}\left(#1 \right)}}
\newcommand{\np}{\ensuremath{\mathcal{NP}}}

\newtheorem{proposition}{Proposition}
\newtheorem{theorem}{Theorem}
\newtheorem{corollary}{Corollary}
\newtheorem{lemma}{Lemma}
\theoremstyle{definition}
\newtheorem{definition}{Definition}

\theoremstyle{remark}
\newtheorem{remark}{Remark}

%% file: notations.tex
\newcommand\defeq{\stackrel{\mathclap{\text{\tiny def}}}{=}}
\DeclareMathOperator{\rank}{rank}

\DeclareMathOperator{\support}{supp}
\newcommand*{\ones}{\ensuremath{\mathbf{e}}}
\newcommand{\littleo}{\ensuremath{\scalebox{.9}{$\scriptscriptstyle\mathcal{O}$}}}
\newcommand{\setRankOne}{\set{Q}_1}